\documentclass[smallextended,referee,envcountsect]{svjour3}
\usepackage[pdftex]{graphicx}
\usepackage{graphicx}
\usepackage{float}
\usepackage{color}
\usepackage{hyperref}
\hypersetup{colorlinks=true,
            linkcolor={blue},
            citecolor={green}}

\usepackage{amsmath, amsfonts, latexsym}
\usepackage{paralist}
\usepackage{mathabx}
\usepackage{mathrsfs,amssymb}
\usepackage{version}
\usepackage{bm}
\usepackage{bbm}
\usepackage{flafter}
\usepackage{ulem}
\usepackage{dsfont}
\usepackage[english]{babel}
\usepackage{cancel}

\newcommand{\cD}{\mathcal{D}}

\newcommand{\cS}{\mathcal{S}}

\newcommand{\cX}{\mathcal{X}}
\newcommand{\cY}{\mathcal{Y}}

\newcommand{\E}{\mathbb{E}}
\newcommand{\F}{\mathbb{F}}

\newcommand{\R}{\mathbb{R}}

\newcommand{\1}{{\bf 1}}
    
\newcommand{\HP}[1] 
    {\ensuremath{\mathscr{H}^{#1}}}

\newcommand{\p}{\mathbb P}
\newcommand{\N}{\mathbb N}

\newcommand{\s}{\mathbb S}
\newcommand{\defeq}{\mathrel{\mathop:}=}

\DeclareMathOperator{\interior}{int}

\DeclareMathOperator{\Tr}{Tr}

 \newtheorem{assump}{Assumption}[section]

 \newcommand{\lf}[1]{ #1^\sharp}

 \newcommand{\pr}[1]{ #1^\flat}

\smartqed

\begin{document}

 \title{A level-set approach for stochastic optimal control problems under controlled-loss constraints}

\author{G\'{e}raldine Bouveret   \and  Athena Picarelli }

\institute{G\'{e}raldine Bouveret \at
             School of Physical and Mathematical Sciences, Nanyang Technological University, Singapore\\
              geraldine.bouveret@ntu.edu.sg
           \and
           Athena Picarelli  \at
              Department of Economics, Universit\`a di Verona, Verona, Italy\\
             athena.picarelli@univr.it
}

\date{\textcolor{white}{Received: date / Accepted: date}}

 \maketitle

 \begin{abstract}
We study a family of  optimal control problems under a set of controlled-loss constraints holding at different deterministic dates. The characterization of the associated value function by a Hamilton-Jacobi-Bellman equation usually calls for additional strong assumptions on the dynamics of the processes involved and the set of constraints. To treat this problem in absence of those assumptions, we first convert it into a state-constrained stochastic target problem and then apply a level-set approach. With this approach, the state constraints can be managed through an exact penalization technique.
    \end{abstract}


  \section{Introduction}
Under some general assumptions, the value function associated with stochastic optimal control problems can be characterized
as the unique continuous  viscosity solution of a second-order Hamilton-Jacobi-Bellman (HJB) equation (see e.g. \cite{YZ99,FS06} and the references therein).
However in many applications, the question of optimization under state constraints arises.
We suggest here to study the case of state constraints holding in expectation and on a set of deterministic dates. Those constraints involve a loss function and are thus referred in the literature  as \textit{controlled-loss constraints}. They are of great interest in finance (see e.g. \cite{ELJBL05,BT07,GW07,DT11}).

There is a huge literature on state-constrained optimal control problems and their associated HJB equations (see e.g. \cite{K94,BR98,IL02,F08,BEI10}).
The characterization of the value function as a viscosity solution of a HJB equation is intricate and usually requires a delicate interplay between the dynamics of the processes involved and the set of constraints. In particular, some viability and regularity assumptions on the dynamics are typically required to ensure the finiteness of the value function and its partial differential equation (PDE) characterization, often rendering the problem not treatable.

From a mathematical perspective, the state-constrained problem considered here is \textit{non-standard} as the constraints are expressed in expectation and imposed at different discrete times. Similar problems have been studied in \cite{B18} with only one probabilistic constraint, and the analysis still involves strong assumptions on the controls and processes involved. In particular, controls must be real-valued.
In the discrete-time setting, optimal control problems with one expectation constraint are studied in \cite{P18}.

This paper therefore aims at providing an alternative way for solving such problems under a more general framework.
This objective is achieved  at the price of augmenting the state and control space by additional components and considering unbounded controls.
More precisely, following the ideas developed in \cite{BPZ16} for the case of a constraint holding pointwise almost surely, our approach relies on two main steps. The first one consists in building on the arguments developed in \cite{BET09,BD12} to reformulate, by means of the martingale representation theorem, the original problem as a stochastic target problem involving almost-sure constraints and unbounded controls.
Stochastic target problems have been extensively studied via a HJB characterization \cite{ST02GEO,ST02,BET09,BC17} or a dual approach \cite{BBC16,B18b}, the latter suffering from a lack of tractability of computations when several dates are involved in the constraint. Still, a direct treatment of the derived stochastic target problem remains challenging due to the nature and number of constraints involved (see e.g. \cite{BEI10,BD12,B18}).
The second step therefore consists in solving  the resulting stochastic target problem by means of a \textit{level-set approach} where the state constraints are managed via an exact penalization technique. Initially introduced in \cite{OS88}
to model some deterministic front  propagation phenomena, the level-set approach has been used in many applications related to (non)linear controlled systems (see e.g. \cite{FGL94,KV06,BFZ10,GP15}). The connection between stochastic target problems and level-set characterization has been already pointed out in \cite{ST02b}.
In our case, the level-set approach links  the stochastic target  problem to an auxiliary optimal control problem, referred as the \textit{level-set problem}, defined on the augmented state and control space, but without state constraints. The value function associated with this level-set problem can be fully characterized as the solution of a particular HJB equation, offering a complete treatment of the original problem. In particular, due to the presence of unbounded controls, a compactification of the differential operator must be applied as in \cite{BPZ16,BC17} to derive comparison results.
The main difference with \cite{BPZ16} relies on the time inconsistency of the problem. Indeed, at a specific time, the number of unbounded controls involved depends on the number of constraints holding in the future and thus on the time interval considered. The associated value function is therefore defined differently on each time interval $[t_i, t_{i+1})$, $0\leq i\leq n-1,$ and shows discontinuity at $t_{i+1}$.
Additionally, this approach relies on a relaxed version of the condition of existence of an optimizer for the level-set problem stated in \cite[(H4)]{BPZ16}, which expands the scope of application of the results. Although this condition relates to convexity properties of the dynamics, cost functions, and the set of controls, it is  independent from  the aforementioned viability and regularity assumptions needed to directly characterize the original problem. Finally, uniform boundedness in $L^2$ of the admissible controls is not considered here and the regularity of the value function associated with the level-set problem cannot be proven a priori.

The rest of the paper is organized as follows. In Section \ref{pbstatement}, we formally state the problem. In Section \ref{se PbRdPbedp}, we formulate the optimal control problem as a constrained stochastic target problem. The level-set approach is then applied provided that an assumption of control existence is satisfied. The latter assumption is investigated in Section \ref{se existence}. A complete characterization of the obtained level-set function is derived in Section \ref{se charac}. An appendix contains proofs of some technical results.

\textbf{Notations.} We let $d,q\ge 1$ be integers. Every element of $\R^d$ is considered as a column vector. We denote by $|x|$ the Euclidean norm of $x\in \R^d$, and by $x^{\top}$ its transpose. We also set $\mathrm M^{\top}$ the transpose of $\mathrm M\in\R^{d\times q}$, while $\Tr[\mathrm M]$ is its trace.   We  denote by  $\s^{d}$  the set
 of symmetric matrices in $\R^{d\times d}$ and by $I_d\in\s^{d}$ (resp. $\textbf{0}\in \s^{d}$), the identity matrix (resp. the null matrix).
Moreover, we define $\mathcal{S}_d$ the unit $d$-sphere, i.e. $\{b\in\R^{d+1},\,|b|=1\}$, and $\mathcal{D}_d$ the subset of $\mathcal{S}_d$ such that the first component $b_1$ is null.
Finally, $\R^+:=[0,\infty)$ and $\R^+_\ast:=(0,\infty)$, $\R^-\defeq(-\infty,0]$, $C, \hat C,\bar C>0$ are constant terms that we do not keep track of, the abbreviation ``s.t.'' stands for ``such that'', and  inequalities between random variables hold $\p$-a.s.

\section{Setting and main assumptions}\label{pbstatement}
   In this manuscript we consider $\Omega$, the space of $\R^q$-valued continuous functions $(\omega_t)_{t\leq T}$ on $[0,T]$ endowed with the Wiener measure $\p$. We introduce $W$ the coordinate mapping, i.e. $(W(\omega)_t)_{t\leq T}$ for $\omega \in \Omega$ so that $W$ is a $q$-dimensional Brownian motion on the canonical filtered probability space $(\Omega,\mathcal{F},\F,\p)$. In particular, $\mathcal{F}$ is the Borel tribe of $\Omega$ and $\F\defeq\{\mathcal{F}_t,0\leq t\leq T\}$ is the $\p$-augmentation of the filtration generated by $W$.
   We define $\mathcal{U}$ as the collection of progressively measurable processes $\nu$ with values in $U$, a compact subset of $\R^r,\,r\ge 1$.
   For $t\in [0,T], z\in \R^{d}$ and for $\nu \in \mathcal{U}$, we define the process $Z_\cdot^{t,z,\nu}$
as the unique (strong) solution on $[t,T]$ to
$$ Z^{t,z,\nu}_\cdot=z+\int_{t}^{\cdot}\mu(s,Z^{t,z,\nu}_s, \nu_s) \,  \mathrm{d}s  +\int_{t}^{\cdot} \sigma(s,Z^{t,z,\nu}_s,\nu_s) \, \mathrm{d}W_s\,,$$
 where $(\mu,\sigma): (t,z,u)\in[0,T]\times\R^d\times U \rightarrow \R^d\times \R^{d\times q}$  are continuous functions being Lipschitz continuous in $z$ and satisfying a linear growth in $z$ uniformly in $(t,\nu)$.
\begin{remark}\label{rem:XY}
Financial applications usually consider the case where $Z^{t,z,\nu}_\cdot \defeq (X^{t,x,\nu}_\cdot, Y^{t,z,\nu}_\cdot)$, with
  $X^{t,x,\nu}_\cdot\defeq x+\int_{t}^{\cdot}\mu_X(s,X^{t,x,\nu}_s,\nu_s) \,  \mathrm{d}s  +\int_{t}^{\cdot} \sigma_X(s,X^{t,x,\nu}_s,\nu_s) \, \mathrm{d}W_s$ on $\R^{d}$ and $Y^{t,z,\nu}_\cdot\defeq y+\int_{t}^{\cdot} \mu_Y(s,Z^{t,z,\nu}_s,\nu_s) \, \mathrm{d}s +\int_{t}^{\cdot} \sigma_Y^\top (s,Z^{t,z,\nu}_s,\nu_s) \, \mathrm{d}W_s$ on $\R$
where  $(\mu_X,\sigma_X): (t,x,u)\in[0,T]\times\R^d\times U \rightarrow \R^d\times\R^{d\times q}$ (resp. $(\mu_Y,\sigma_Y): (t,z,u)\in[0,T]\times\R^{d+1}\times U \rightarrow \R\times\R^q$) are continuous functions being Lipschitz continuous in $x$ (resp. $z$).
In this form, $X$ models the evolution over time of the price of some underlying assets while $Y$ represents a portfolio process.
\end{remark}
To simplify the notations, we assume from now on that the dimensions $d$, $q$ and $r$ are all equal and we thus disregard the notations $q$ and $r$.

We now introduce two non-negative Lipschitz continuous maps  $f$ and $\Psi$ defined on $\R^d$. We fix $n\in \N$ and consider the time grid $t_0=0\leq \dots \leq t_i\leq \dots\leq t_n=T$.

For any $0\le i\le n-1$, we define the set $\mathscr{C}_i\defeq[t_i,t_{i+1})\times\R^d\times \R^{(n-i)},$ as well as
 \begin{align*}
& \mathscr{B}_i\defeq[t_i,t_{i+1})\times\R^{d}\times[0,\infty)^{n-i}, && \interior(\mathscr{B}_i)\defeq[t_i,t_{i+1})\times\R^{d}\times(0,\infty)^{n-i}\,,\\
 &\mathscr{D}_i\defeq \mathscr{B}_i\times\R^+, && \interior(\mathscr{D}_i)\defeq\interior(\mathscr{B}_i)\times\R^+_*\,.
\end{align*}
    The objective of the paper is to solve for any $0\le i\le n-1$ the following stochastic optimal control problem on $\mathscr{C}_i$,
\begin{equation}\label{pbintialedp}
  V(t,z,p_{i+1},...,p_n)\defeq\inf_{\nu\in{\mathcal{U}}_{t,z,p_{i+1},...,p_n}}\E\left[f(Z^{t,z,\nu}_T)\right]\,,
\end{equation}
where
${\mathcal{U}}_{t,z,p_{i+1},...,p_n}\defeq\left\{\nu\in\mathcal{U}:\;\E\left[\Psi( Z^{t,z,\nu}_{t_k})\right]\le p_k,\,i+1\le k\le n\right\}.
$
On $\{T\}\times\R^{d}$, we set $V(T,z)=f(z)$. We use the convention $V(t,z,p_{i+1},...,p_n)=\infty$ whenever ${\mathcal{U}}_{t,z,p_{i+1},...,p_n}=\emptyset$.
Observe that ${\mathcal{U}}_{t,z,p_{i+1},...,p_n}=\emptyset$ whenever there exists $i+1\le k\le n$ s.t. $p_k<0$. This implies that $V=\infty$ on $\mathscr{C}_i\setminus \mathscr{B}_i$. We underline that the problem can be treated similarly if we consider different loss functions at each date.
\begin{remark}
Consider some continuous non-negative function $\ell$ defined on $[0,T]\times\R^{d}\times U$ and being Lipschitz continuous in the space variable. One can always consider  the augmented dynamics  $\tilde Z^{t,z,\nu}_{\cdot} \defeq (Z^{t,z,\nu}_\cdot, \zeta^{t,z,\nu}_\cdot)\in \R^{d+1}$ for $\zeta^{t,z,\nu}_\cdot\defeq \int^\cdot_t \ell (s,Z^{t,z,\nu}_s,\nu_s) \mathrm d s$, together with the terminal condition $\tilde f(\tilde Z^{t,z,\nu}_T) \defeq f(Z^{t,z,\nu}_T)+\zeta^{t,z,\nu}_T$, and recover the formulation of the problem given in \eqref{pbintialedp}.
\end{remark}
 \section{Problem reformulation}\label{se PbRdPbedp}
In the spirit of \cite{BPZ16}, our approach articulates in two steps. First, we reformulate  \eqref{pbintialedp} as a constrained stochastic target problem (see Proposition \ref{prop stotget} below).  Then, this  stochastic target problem is described
by a level-set approach where the constraints are handled using an exact penalization technique (see Proposition \ref{prop aux} below).
This links the backward reachable set associated with the stochastic target problem to the zero level-set of a value function associated with a suitable auxiliary unconstrained optimal control problem given by $w$ in \eqref{defw} below.
\subsection{Associated stochastic target problem}\label{se target}
We denote by  $\mathcal{A}$, the collection of progressively measurable processes in $L^2([0,T]\times\Omega)$, with values in $\R^d$.
Let $0\leq i\leq n-1$ and $i+1\le k\le n$. Before stating the main result, we define for any  $t\in [t_i,t_{i+1})$, $p_k\in\R$, $\alpha_k\in \mathcal{A}$, $m\in \R$, and $\eta\in \mathcal A$ the new processes
\begin{align*}
P^{t,p_k,\alpha_k}_{\cdot}\defeq p_k+\int_t^{\cdot}\alpha^{\top}_{k_s}\,\mathrm{d}W_s \mbox{ on }[t,t_k]\quad\mbox{ and }\quad M^{t,m,\eta}_{\cdot}\defeq m+\int_t^{\cdot}\eta_{s}^\top\,\mathrm{d}W_s \mbox{ on }[t,T]\,.
\end{align*}
We shortly denote $\alpha\equiv(\alpha_{i+1},\ldots,\alpha_n)\in \mathcal A\times \ldots\times\mathcal A \equiv \mathcal A^{n-i}$.
We can now state the following result.
\begin{proposition}\label{prop stotget}
Let  $0\le i\le n-1.$ For any $(t,z,p_{i+1},...,p_n)\in\mathscr{C}_i,$
\begin{align}
 V(t,z,p_{i+1},...,p_n) =  \inf\bigg\{ & m\ge 0:\exists\,(\nu,\alpha,\eta)\in {\mathcal{U}}\times{\mathcal{A}}^{n-i}\times\mathcal{A}\mbox{ s.t. } \label{eq defV}\\
& \hspace{-0.5 cm}M^{t,m,\eta}_T\ge f(Z^{t,z,\nu}_T)\; \text{and}\;
P^{t,p_k,\alpha_k}_{t_k}\ge \Psi(Z^{t,z,\nu}_{t_k} ),\,i+1\le k\le n \nonumber
\bigg\}\,.
\end{align}
\end{proposition}
\noindent{\it Proof}
Let $0\le i\le n-1$ and $(t,z,p_{i+1},...,p_n)\in\mathscr{C}_i$.
One can easily prove that
\begin{align*}
V(t,z,p_{i+1},...,p_n)=\inf\left\{m\ge 0: \exists\,\nu\in\mathcal{U}_{t,z,p_{i+1},...,p_n}\mbox{ s.t. }m\ge \E\left[f(Z^{t,z,\nu}_T)\right]\right\}\,.
\end{align*}
We then prove for any $m\ge 0$, the equivalence between the two following statements
\begin{align*}
    (i)& \,\exists \,\nu\in\mathcal{U}_{t,z,p_{i+1},...,p_n}\mbox{ s.t. }m\ge \E\left[f(Z^{t,z,\nu}_T)\right] \,,\\
    (ii)&\,\exists \,(\nu,\alpha,\eta)\in{\mathcal{U}}\times\mathcal{A}^{n-i}\times\mathcal{A},\mbox{ s.t. }\begin{cases} P^{t,p_k,\alpha_k}_{t_k}\ge \Psi(Z^{t,z,\nu}_{t_k}),\,i+1\le k\le n\\\mbox{and }M^{t,m,\eta}_T\ge f(Z^{t,z,\nu}_T)\end{cases}\,.
\end{align*}
To this aim we appeal to similar techniques as those exploited in \cite{BET09,BBC16}.
The implication $(ii)\Rightarrow (i)$ follows by taking the expectation in $(ii)$ and using the martingale property of the stochastic integrals. On the other hand, the implication $(i)\Rightarrow (ii)$ follows from the martingale representation theorem (see e.g. \cite[Theorem 4.15, Chapter 3]{KS12}). More precisely, from the assumptions on the coefficients of $Z$ and the growth conditions on $f$ and $\Psi$, there exists, for any $\nu\in\mathcal{U}$, $(\hat{\alpha}_{i+1},...,\hat{\alpha}_{n},\hat\eta)\in\mathcal{A}^{n-i}\times\mathcal{A}$ s.t.
\begin{align*}
M^{t,m,\hat\eta}_T= m+\int^T_{t}\hat\eta^\top_s\,\mathrm{d} W_s\ge f(Z^{t,z,\nu}_T)
\mbox{ and }
P^{t,p_k,\hat\alpha_k}_{t_k}=p_k+\int^{t_k}_{t}\hat\alpha^{\top}_{k_s}\,\mathrm{d} W_s\ge \Psi(Z^{t,z,\nu}_{t_k}),\,i+1\le k\le n
\,,
\end{align*}
leading to the result.
\qed

As intimated  in the introduction, a direct treatment of the derived stochastic target problem \eqref{eq defV} is challenging and would involve strong regularity assumptions that would considerably restrict the applicability of our study (see e.g. \cite{BEI10,BD12,B18}). Accordingly, a direct resolution of the derived stochastic target problem seems unsatisfactory and we will thus link, in Section \ref{se compact}, $V$ to the auxiliary problem given by $w$ in \eqref{defw} below.

\subsection{Level-set approach}\label{se compact}
For $0\le i\le n-1,$ and $(t,z,p_{i+1},...,p_n,m)\in \mathscr{C}_i\times\R$, we define the following optimal control problem
\begin{gather}
w(t,z,p_{i+1},...,p_n,m)\defeq
\inf_{\substack{{(\nu,\alpha,\eta)\in\mathcal{U}\times {\mathcal{A}}^{n-i}\times\mathcal{A}}}}{\mathrm{J}}^{{\nu},\alpha,\eta}(t,z,p_{i+1},...,p_n,m)
\,,\label{defw}
\end{gather}
with
\begin{align*}
{\mathrm{J}}^{{\nu},\alpha,\eta}(t,z,p_{i+1},...,p_n,m)\defeq\E\left[\left(f(Z^{t,z,\nu}_{T})-M^{t,m,\eta}_T\right)^+
+\sum^{n}_{k=i+1}\left(\Psi(Z^{t,z,\nu}_{t_{k}})-P^{t,p_{k},\alpha_{k}}_{t_{k}}\right)^+\right]\,.\end{align*}
On $\{T\}\times\R^{d}\times\R$, we set $w(T,z,m)=(f(z)-m)^+$.

In what follows, we denote $w^*$ (resp. $w_*$) the upper (resp. lower) semi-continuous envelope of $w$ on $\mathscr{D}_i,\,0\le i\le n-1$.
The function $w$ satisfies the following regularity properties.
\begin{proposition}\label{prop growth}
For any $0\le i\le n-1$, $w$ is Lipschitz continuous with respect to $(z,p_{i+1},...,p_n,m)$ on $\mathscr C_i\times\R$ and satisfies on $\mathscr D_i$,
\begin{align}\label{eq polygrowth}
0\leq w(t,z,p_{i+1},...,p_n,m)\le C\left(1+|z|\right)\,,
\end{align}
for some $C>0$.
Moreover,  on $\R^d\times\R\times\R$,
\begin{align}\label{eq termcondwtilde}
\lim_{t\uparrow T}w(t,z,p_n,m)=  (f(z)-m)^+ + (\Psi(z)-p_{n})^+\,.
\end{align}
\end{proposition}
\noindent{\it Proof} Fix $0\le i\le n-1$. The Lipschitz continuity of $w$ with respect to $(p_{i+1},...,p_n,m)$ (resp. $z$) is straightforward (resp. follows from the regularity of $f$, $\Psi$ and of the coefficients of $Z$).
 Moreover, by definition of $w$, one has on $\mathscr{C}_i\times\R$,
\begin{gather*}
0\leq w(t,z,p_{i+1},...,p_n,m)\le \inf_{\substack{{\nu\in\mathcal{U}}}}\,\E\left[\left(f(Z^{t,z,\nu}_T)-m\right)^+
+\sum_{k=i+1}^n\left(\Psi(Z^{t,z,\nu}_{t_k})-p_k\right)^+\right]\,.
\end{gather*}
Therefore, since $m,p_k\geq 0,\,i+1\le k\le n,$ on $\mathscr D_i $,
\eqref{eq polygrowth} follows from the growth conditions on $f$ and $\Psi$ and the assumptions on the coefficients of $Z$.
We now prove \eqref{eq termcondwtilde}.  Let $0< h\le T-t_{n-1}$.
On $\R^d\times\R\times\R$, one has for any $\nu\in \mathcal U$,
\begin{align*}
& {w}(T-h,z,p_{n},m) \leq \E\left[\left(f(Z^{T-h,z,\nu}_T)-m\right)^++\left(\Psi(Z^{T-h,z,\nu}_T)-p_n\right)^+\right]
\,.
\end{align*}
On the other hand, by the martingale property of stochastic integrals, one has
\begin{align*}
 {w}(T-h,z,p_{n},m)\geq  \left(\inf_{\nu\in \mathcal U} \E\left[ f(Z^{T-h,z,\nu}_T)\right]-m\right)^++\left(\inf_{\nu\in \mathcal U} \E\left[\Psi(Z^{T-h,z,\nu}_{T}) \right]-p_n\right)^+ \,.
\end{align*}
Therefore, the Lipschitz continuity of  $f$ and $\Psi$ together with  classical estimates on the process $Z$, provide the existence of a uniform $C>0$ s.t.
\begin{align}\label{eq lim2}
|{w}(T-h,z,p_{n},m) -  (f(z)-m)^+ - (\Psi(z)-p_{n})^+|\le  C\sqrt h \left(1+|z|\right)\,.
\end{align}
We finally let $h$ tend to zero in \eqref{eq lim2} to conclude.
\qed

The following assumption is a  weaker condition than \cite[(H4)]{BPZ16} and is key for proving Proposition \ref{prop aux} below.
\begin{assump}\label{ass ctrolexistence}
On $\mathscr{C}_i\times\R,\,0\le i\le n-1$, if $w(t,z,p_{i+1},\ldots,p_n,m)=0$ at some point  $(t,z,p_{i+1},\ldots, p_n, m)$, then there exists an optimal control for the problem \eqref{defw} at  $(t,z,p_{i+1},\ldots, p_n, m)$.
\end{assump}

\begin{remark}\label{re boundctrol}
Assumption \ref{ass ctrolexistence} is weaker than \cite[(H4)]{BPZ16} since it only requires the existence of an optimal control at those points  $(t,z,p_{i+1},\ldots, p_n, m)$ where $w$ is zero. This is the minimal requirement to obtain the characterization of $V$ in Proposition \ref{prop aux} below. Differently from \cite[(H4)]{BPZ16},  Assumption \ref{ass ctrolexistence} can be proved to hold, under suitable convexity assumption (see next section), without  requiring a uniform bound in the $L^2$-norm of controls. 
\end{remark}

\begin{proposition}\label{prop aux}
Under Assumption \ref{ass ctrolexistence}, one has  on $\mathscr{C}_i,\,0\le i\le n-1$,
\begin{align*}
V(t,z,p_{i+1},...,p_n)=\inf\left\{m\ge 0:  w(t,z,p_{i+1},...,p_n,m)=0\right\}\,.
\end{align*}
\end{proposition}
\noindent{\it Proof}
Fix $0\le i\le n-1$. In virtue of Proposition \ref{prop stotget} it is sufficient to show that for any point $(t,z,p_{i+1},...,p_n,m)\in \mathscr{C}_i\times \R^+$ the following equivalence holds
\begin{align}\label{eq Identity}
& \exists (\nu,\alpha,\eta)\in\mathcal{U}\times\mathcal{A}^{n-i}\times \mathcal A\mbox{ s.t. } M^{t,m,\eta}_T\ge f(Z^{t,z,\nu}_T)\mbox{ and }
P^{t,p_k,\alpha_k}_{t_k}\ge \Psi(Z^{t,z,\nu}_{t_k}),\quad i+1\le k\le n\nonumber\\
& \Leftrightarrow w(t,z,p_{i+1},...,p_n,m)=0\,.
\end{align}
\textbf{Step 1. Proof of $\Rightarrow$.} The implication follows after observing that when the left-hand side of the equivalence in \eqref{eq Identity} is satisfied,
one has
\begin{gather}
\left(f(Z^{t,z,\nu}_T)-M^{t,m,\eta}_T\right)^+=0\mbox{ and }
\left(\Psi(Z^{t,z,\nu}_{t_{k}})-P^{t,p_{k},\alpha_{k}}_{t_{k}}\right)^+=0,\quad i+1\le k\le n\,,\label{positivity2}
\end{gather}
leading to $w(t,z,p_{i+1},...,p_n,m)=0$.\\
\textbf{Step 2. Proof of $\Leftarrow$.}
Let $ w(t,z,p_{i+1},...,p_n,m)=0$. We appeal to Assumption \ref{ass ctrolexistence} and consider the  optimal control $(\nu,\alpha,\eta)\in\mathcal{U}\times\mathcal{A}^{n-i}\times \mathcal A$ s.t.
$$
\E\left[\left(f(Z^{t,z,\nu}_{T})-M^{t,m,\eta}_T\right)^+
+\sum^{n}_{k=i+1}\left(\Psi(Z^{t,z,\nu}_{t_{k}})-P^{t,p_{k},\alpha_{k}}_{t_{k}}\right)^+\right]=0\,.
$$
 Thus  \eqref{positivity2} holds by the non-negativity of each term, hence the implication.
\qed
\begin{remark}\textcolor{white}{}
One can easily verify that  $\inf\left\{m\ge 0:  w(T,z,m)=0\right\}=f(z)=V(T,z)$ and thus that the result in Proposition \ref{prop aux} extends to $\{T\}\times\R^{d}$.
\end{remark}
Proposition \ref{prop aux} is critical here as it allows the reformulation of $V$ in terms of the unconstrained optimal control problem described by $w$ whose associated value function satisfies important regularity properties (recall Proposition \ref{prop growth}). Therefore a complete PDE characterization can be provided in Section \ref{se charac}. We point out that problem \eqref{defw} is a singular optimal control problem characterized by a discontinuous Hamiltonian. As a result, the HJB characterization must be obtained passing through a reformulation of the differential operator for a comparison result to hold (see e.g. \cite{BBMZ09,BPZ16,BC17}). Observe also that, unlike \cite{BPZ16}, the cost functional associated with $w$ changes on each time interval $[t_i, t_{i+1}),\,0\le i\le n-1,$ to adapt to the decreasing number of constraints involved. As a result, a discontinuity at each point $t_i,\,1\le i\le n,$ arises.

\section{Existence results}\label{se existence}
We give in this section some  sufficient conditions ensuring that Assumption  \ref{ass ctrolexistence} is satisfied (recall Remark \ref{re boundctrol}).
\begin{proposition}\label{prop:exists}
Assume that $U$ is a convex set, $f$ and $\Psi$ are convex functions, and the coefficients of the diffusion are of the form
$
\mu(t,z,u)\defeq A(t)z+B(t)u$ and $\sigma(t,z,u)\defeq C(t)z+D(t)u,
$
for all $(t,z,u)\in[0,T]\times\R^d\times U$, with $A, B, C$ and $D$ matrices of suitable size. If  $(t,z,p_{i+1},\ldots,p_n,m) \in \mathscr C_i\times \R$, for some $0\leq i\leq {n-1}$ is  such that $w(t,z,p_{i+1},\ldots,p_n,m)$ $=0$, then, the optimal control problem \eqref{defw} admits an optimizer at $(t,z,p_{i+1},\ldots,p_n,m)$.
\end{proposition}
\noindent{\it Proof}
Fix $0\le i\le n-1$. Let $(\nu^j,\alpha^j,\eta^j)\in \mathcal U\times \mathcal A^{n-i} \times \mathcal A$ be a minimizing sequence for  $w$ at point $(t,z,p_{i+1},\ldots,p_n,m)\in\mathscr{C}_i\times\R$ s.t. $w(t,z,p_{i+1},\ldots,p_n,m) = 0$.
Therefore, for any $\varepsilon>0$  there exists $j_0$ s.t. $\forall j\geq j_0$ one has
$$
\mathbb E\bigg[\big(f(Z^{t,z,\nu^j}_T) -  M^{t,m,\eta^j}_T\big)^+ + \sum^n_{k=i+1} \big(\Psi(Z^{t,z,\nu^j}_{t_k} ) -  P^{t,p_k,\alpha^j_k}_{t_k}\big)^+ \bigg] \leq w(t,z,p_{i+1},\ldots,p_n,m) + \varepsilon\,.
$$
As $\nu^j$ is uniformly bounded in the $L^2$-norm (since $\nu$ takes values in the compact set $U$), there exists a subsequence (still indexed by $j$)  $\nu^j$ that weakly converges in the $L^2$-norm to some $\hat \nu\in \mathcal U$.
Applying Mazur's theorem, one has the existence of $\tilde \nu^j\equiv \sum_{\ell\geq 0} \lambda_\ell \nu^{\ell+j}$ with $\lambda_\ell\geq 0$ and $\sum_{\ell\geq 0}\lambda_\ell =1$ s.t. $\tilde \nu^j$ strongly converges in the $L^2$-norm to $\hat \nu$.
We then consider $(\tilde \nu^j, \tilde \alpha^j, \tilde\eta^j)\equiv \sum_{\ell\geq 0} \lambda_\ell (\nu^{\ell+j},\alpha^{\ell+j}, \eta^{\ell+j})$. Observe that $(\tilde \nu^j, \tilde \alpha^j, \tilde \eta^j)$ still belongs to $\mathcal U\times  \mathcal A^{n-i}\times \mathcal A$ as $\mathcal{U}$ and $\mathcal A$ are convex spaces.  \\
 Let us now consider $(Z^{t,z,\tilde \nu^j}_\cdot, M^{t,m,\tilde \eta^j}_\cdot, (P^{t,p_k,\tilde \alpha^j_k}_\cdot)_{i+1\leq k\leq n})$. One has for $j\to+\infty$,
  $$
  \mathbb E\left[\sup_{s\in [t,T]} \left| Z^{t,z,\tilde \nu^j}_s - Z^{t,z,\hat \nu}_s\right|^2\right] \to 0\,.
  $$
It follows that  for any $\varepsilon >0$ there exists $j_1$ s.t.  for all $j\geq j_1$,
$$
\left | \mathrm{J}^{\tilde \nu^j,\tilde \alpha^j, \tilde \eta^j}(t,z,p_{i+1},...,p_n,m) -  \mathrm{J}^{\hat \nu,\tilde \alpha^j, \tilde \eta^j}(t,z,p_{i+1},...,p_n,m)\right| \leq \varepsilon\,.
$$
 Moreover, by the linearity of the dynamics of $Z$ one has
$
 Z^{t,z,\tilde \nu^j}_\cdot  = \sum_{\ell\geq 0} \lambda_\ell Z^{t,z,\nu^{\ell+j}}_\cdot\,,
$
and trivially
  $
 M^{t,m,\tilde \eta^j}_\cdot = \sum_{\ell\geq 0} \lambda_\ell  M^{t,m,\eta^{\ell+j}}_\cdot \quad\text{and}\quad  P^{t,p_k,\tilde \alpha^j_k}_\cdot = \sum_{\ell\geq 0} \lambda_\ell  P^{t,p_k,\alpha^{\ell+j}_k}_\cdot \text{ for } i+1\leq k\leq n\,.
 $
 Therefore, using the convexity of $f$ and $\Psi$, one has for $j\geq \max(j_0,j_1)$,
 \begin{align*}
&  \mathbb E \left[\big( f(Z^{t,z,\hat \nu}_T) - M^{t,m,\tilde \eta^j}_T \big)^+  + \sum^n_{k=i+1} \big(  \Psi(Z^{t,z,\hat \nu}_{t_k}) - P^{t, p_k,\tilde \alpha^j}_{t_k} \big)^+ \right] \\
 &\quad \quad\leq \sum_{\ell\geq 0} \lambda_{\ell} \mathbb E\bigg[\big( f(Z^{t,z, \nu^{\ell+j}}_T) - M^{t,m,\eta^{\ell+j}}_T \big)^+ +\sum^n_{k=i+1} \big(  \Psi(Z^{t,z,\nu^{\ell+j}}_{t_k}) - P^{t, p_k, \alpha^{\ell+j}}_{t_k} \big)^+\bigg] + \varepsilon\\
 & \quad\quad\leq w(t,z,p_{i+1},\ldots,p_n,m) +2 \varepsilon
=2\varepsilon\,.
 \end{align*}
 Thanks to the arbitrariness of $\varepsilon$ and  the martingale property of stochastic integrals, it is immediate to verify that the previous inequality gives
\begin{equation}\label{exist_ineq}
\mathbb E\left[ f(Z^{t,z,\hat \nu}_{T} ) \right]\leq m \qquad\text{and}\qquad \mathbb E\left[ \Psi(Z^{t,z,\hat \nu}_{t_k} ) \right]\leq p_k \quad \text{for any } i+1\leq k\leq n\,.
\end{equation}
By the martingale representation theorem, we define  $\hat \eta,\hat \alpha_k \in \mathcal A,\,i+1\leq k\leq n$, s.t.
$$
\mathbb E\left[ f(Z^{t,z,\hat \nu}_{T} ) \right] = f(Z^{t,z,\hat \nu}_{T} ) - \int^T_t \hat \eta_s^\top \mathrm d W_s \;\;\;\text{and}\;\;\; \mathbb E\left[ \Psi(Z^{t,z,\hat \nu}_{t_k} ) \right] = \Psi(Z^{t,z,\hat \nu}_{t_k} ) - \int^{t_k}_t \hat \alpha_{k_s}^\top \mathrm d W_s\,.
$$
In virtue of \eqref{exist_ineq}, this gives
 \begin{align*}
\mathrm{J}^{\hat \nu,\hat \alpha, \hat \eta}(t,z,p_{i+1},...,p_n,m)   =   \mathbb E \left[\big( f(Z^{t,z,\hat \nu}_T) - M^{t,m,\hat \eta}_T \big)^+  + \sum^n_{k=i+1} \big(  \Psi(Z^{t,z,\hat \nu}_{t_k}) - P^{t, p_k,\hat \alpha_k}_{t_k} \big)^+ \right] = 0\,,
 \end{align*}
from which the optimality of the control $(\hat \nu, \hat \alpha_{i+1},\ldots, \hat\alpha_n,\hat \eta)$ follows.
\qed

For a dynamics $Z\defeq (X,Y)$ (recall Remark \ref{rem:XY}), the result of Proposition \ref{prop:exists} also holds true if
$X$ is independent of $\nu$ and the coefficients of $Y$ are of the form
$$\mu_Y(t,z,u)\defeq A(t,x)y+B(t,x)u\quad\mbox{and}\quad \sigma_Y(t,z,u)\defeq C(t,x)y+D(t,x)u\,,$$
with $A, B, C$ and $D$ matrices of suitable size.

We point out that  Assumption \ref{ass ctrolexistence} is the unique important restriction in our approach. However, such requirement is only related to the convexity properties of the dynamics and cost functions defining the original optimal control problem and does not involve any viability assumption usually necessary to deal with state constrained problems.
\section{A complete PDE characterization for $w$}\label{se charac}
In this section, we characterize $w$ as the viscosity solution of a suitable HJB equation with specific boundary conditions.
We restrict the characterization to $\cup_{0\le i\le n-1}\mathscr{D}_i$, as outside this set $V=\infty$.
\subsection{On the interior of the domain}\label{se interiorHJB}
 The main ingredient towards the PDE characterization of $ w$ is the DPP stated below (see Appendix A for the proof).
\begin{theorem}[DPP]\label{th dpp}
Fix $(t,z,p_{i+1},...,p_n,m)\in\mathscr{C}_i\times\R,\,0\le i\le n-1,$ and let $t_i\le \theta <t_{i+1}$  be a stopping time. Then
\begin{align}\label{eq dpp}
w(t,z,p_{i+1},\hdots,p_n,m)= \inf_{(\nu,\alpha,\eta)\in \mathcal{U}\times\mathcal A^{n-i}\times \mathcal A} \E\left[ w(\theta, {Z}^{t,z,{\nu}}_{\theta}, {P}^{t,p_{i+1},{\alpha}_{i+1}}_{\theta},\hdots ,{P}^{t,p_{n},{\alpha}_{n}}_{\theta},
{M}^{t,m,\eta}_{\theta})\right]\,.
\end{align}
\end{theorem}
We now consider two functions $\kappa:\R\mapsto\R^+_*$ and $\lambda:\R\mapsto\R^+_*$. For  any  $u\in  U$,  $a:= (a_{i+1}, \ldots, a_n)  \in \R^d\times \ldots \times \R^d$, $e\in \R^d$  and for any
$\Theta\defeq(t,z,p_{i+1},...,p_n,m,q,A)\in\mathscr{C}_i\times \R \times \R^{d+n-i+1}\times\s^{d+n-i+1}$, $0\le i\le n-1$, with $q\defeq\left({q^{z\,\top}},q^{p_{i+1}},...,q^{p_n},q^m\right)^\top$, for $q^z\in \R^d$, $q^{p_{k}}\in \R$ $(i+1\leq k \leq n)$, $q^m\in \R$,  and
\begin{align*}
 A & \defeq\left(\begin{matrix}A^{zz}&A^{zp}& A^{zm}\\A^{zp^\top}&A^{pp}& A^{pm}\\A^{zm^\top}&A^{pm^\top}& A^{mm}\end{matrix}\right)\in \s^{d+n-i+1}\,,
 \end{align*}
for $A^{zz}\in \s^d$, $A^{pp}\in \s^{n-i}$, $A^{mm}\in \R$, $A^{zp}\in \R^{d\times (n-i)}$, $A^{zm}\in \R^{d} $ and $A^{pm}\in \R^{n-i}$,
we define  the operators
 \begin{align*}
\mathrm L^{u,a,e}_{\kappa,\lambda} (\Theta)&\defeq
   - \mu^\top(t,z,u)q^z-\frac{1}{2}\Tr[\sigma\sigma^\top(t,z,u)A^{zz}]
  \\
 &\quad -\lambda(m){e^\top}\sigma^\top(t,z,u)A^{zm}  -\sum^{n}_{k=i+1}\kappa(p_k){a^\top_k} \sigma^\top(t,z,u)A^{zp_k}\,,
\\
 \mathrm F^{a,e}_{\kappa,\lambda} (\Theta)&\defeq
 -\frac{1}{2}\lambda(m)^2|e|^2A^{mm} -\frac{1}{2}\sum^{n}_{k=i+1}\kappa(p_k)^2|a_k|^2A^{p_kp_k}
 -\lambda(m)\sum^{n}_{k=i+1}\kappa(p_k){e^\top}{a_k} A^{p_km}\,.
\end{align*}
Hereinafter, we identify each component $a_k\in \R^d,\,i+1\le k\le n,$ of the $(n-i)$-tuple $a$ defined above with the corresponding $d$-dimensional component of the associated vector in $\R^{d(n-i)}$. To alleviate notations we still denote by $a_k,\,i+1\le k\le n,$ (resp. $a$) this component (resp. vector).
 Moreover, for any $\Theta\in\mathscr{C}_i\times \R \times\R^{d+n-i+1}\times\s^{d+n-i+1}$, $c\in \R$, $u\in U$, $b\in \mathcal S_{d(n-i)+d}$, with $b \defeq (b_1, b^{\flat\,\top}_{i+1}, \ldots , {b}^{\flat\,\top}_n , b^{\sharp\,\top})^\top $ for $b_1\in \R$,  $b^\sharp\in \R^d$ and $b^{\flat} := (b^{\flat\,\top}_{i+1} , \ldots , b^{\flat\,\top}_{n})^\top \in \R^{d(n-i)}$,  we also introduce the following operator
 \begin{align*}
 \mathrm H^{u,b}_{\kappa,\lambda}(\Theta,c) \defeq
 \begin{cases}
 (b_1)^2 \left( -c +\mathrm{L}^{u,\bar{b}^{\flat},\bar{b}^{\sharp}}_{\kappa,\lambda}(\Theta) +\mathrm{F}^{\bar{b}^{\flat},\bar{b}^{\sharp}}_{\kappa,\lambda}(\Theta)\right) & b\in \mathcal S_{d (n-i)+d} \setminus \mathcal D_{d (n-i)+d}\\
 \mathrm{F}^{{b}^{\flat},{b}^{\sharp}}_{\kappa,\lambda}(\Theta) & b\in \mathcal D_{d(n-i)+d}
 \end{cases}\,,
 \end{align*}
 where $\bar{b}^\flat\defeq\frac{{b}^\flat}{b_1}=\frac{1}{b_1}(b^{\flat\,\top}_{i+1} , \ldots , b^{\flat\,\top}_{n})^\top\in\R^{d (n-i)}$, $\bar{b}^\sharp\defeq\frac{{b}^\sharp}{b_1}\in\R^d$.\\
Whenever  $\kappa(p)=\lambda(m)=1$ for all $m, p\in \R$, we shortly write $\mathrm F^{a,e}_{\kappa,\lambda} \equiv \mathrm F^{a,e}$, $\mathrm L^{u,a,e}_{\kappa,\lambda} \equiv \mathrm L^{u,a,e}$ and  $\mathrm H^{u,b}_{\kappa,\lambda} \equiv \mathrm H^{u,b}$.
\begin{remark}\label{re opcont}
 The operator $b \mapsto \mathrm{H}_{\kappa,\lambda}^{u,b}$ is continuous
 on $\mathcal{S}_{d (n-i)+d}$, and
 \begin{align*}
 \sup_{\substack{{(u,b)\in U\times\mathcal{S}_{d(n-i)+d}}}} \mathrm H_{\kappa,\lambda}^{u,b}(\Theta,c) = \sup_{\substack{(u,b)\in U\times\mathcal{S}_{d(n-i)+d}\setminus\mathcal{D}_{d(n-i)+d}}} \mathrm H_{\kappa,\lambda}^{u,b}(\Theta,c)\,.
 \end{align*}
 \end{remark}
In what follows, given a smooth function $\varphi$ defined on $\mathscr{C}_i\times \R,\,0\le i\le n-1,$ the notation
$\mathrm H_{\kappa,\lambda}^{u,b}\varphi(\cdot)$  stands for $\mathrm H_{\kappa,\lambda}^{u,b}(\cdot,\mathrm D\varphi(\cdot),\mathrm D^2\varphi(\cdot),\partial_t\varphi(\cdot)).$
A similar writing holds for the operators ${\mathrm L}$ and ${\mathrm F}$.
Finally, observing that the operator $\sup_{(u,a,e)\in  U\times \R^{d(n-i)}\times \R^d}{(\mathrm L^{u,a,e} +\mathrm F^{a,e} )}$ is only lower semi-continuous we denote by
$\big[\sup_{(u,a,e)\in  U\times \R^{d (n-i)}\times \R^d}{(\mathrm L^{u,a,e} +\mathrm F^{a,e} )}\big]^*$ its upper semi-continuous envelope.
 \begin{theorem}\label{theo caractinside}
Let  $0\le i\le n-1$ and define $\kappa,\lambda:\R\mapsto\R^+_*$.  Then, on  $\interior(\mathscr{D}_i)$ the function $w^*$ (resp. $w_*$) is a viscosity sub-solution (resp. super-solution) of
\begin{align}\label{viscsolBETbis}
& \sup_{\substack{{(u,b)\in U\times\mathcal{S}_{d (n-i)+d}}}}\mathrm H_{\kappa,\lambda}^{u,b}w=0\ .
\end{align}
\end{theorem}
\noindent{\it Proof}
Fix $0\le i\le n-1$ and $(t,z,p_{i+1},\ldots,p_n,m)\in \interior(\mathscr{D}_i)$.\\
\textbf{Step 1.}
We first prove the result for $\kappa, \lambda$ identically equal to 1.\\
 We start with the super-solution property.
Let $\varphi$ be a smooth function s.t.
\begin{align}\label{eq:testsuper}
(\text{strict})\min_{\interior(\mathscr{D}_i)}(w_*-\varphi)=(w_*-\varphi)( t, z,p_{i+1},..., p_n,m)=0\,.
\end{align}
 Thanks to Theorem \ref{th dpp}, it follows by standard arguments  (see e.g. \cite[Section 6.2]{BEI10})  that at point  $(t,z, p_{i+1},..., p_n,m)$ one has
\begin{align}\label{supersolBET}
&-\partial_t\varphi+\big[\sup_{(u,a,e)\in  U\times \R^{d (n-i)} \times \R^d}{(\mathrm L^{u,a,e} +\mathrm F^{a,e} )}\big]^*\varphi \geq 0\,.
\end{align}
We then verify that at $( t, z,p_{i+1},..., p_n,m)$,
   $$\sup_{\substack{{(u,b)\in U\times\mathcal{S}_{d(n-i)+d}}}}\mathrm H^{u,b}\varphi \ge 0\,.$$
We adapt the arguments in the proof of \cite[Theorem 3.1]{BC17}.  According to \eqref{supersolBET} and by definition of the upper semi-continuous envelope, we can find a sequence
 $(t_j,z_j,p_{i+1_j},...,p_{n_j},m_j)\in\interior(\mathscr{D}_i),\,q_j\in \R^{d + (n-i) +1}$, $A_j\in\s^{d+n-i+1}$ s.t.
\begin{gather}
(t_j,z_j,p_{i+1_j},...,p_{n_j},m_j)\rightarrow (t, z, p_{i+1},..., p_n, m)\mbox{ and }
|(q_j,A_j)-(\mathrm D\varphi,\mathrm D^2\varphi)(t, z, p_{i+1},...,p_n, m)|\le j^{-1}\,,\label{eq condcv}
\end{gather}
 and
 \begin{align*}
\left(-\partial_t\varphi+\sup_{(u,a,e)\in  U\times \R^{d (n-i) }\times \R^d}\big[\mathrm L^{u,a,e} +\mathrm F^{a,e} \big]\right)(t_j,z_j,p_{i+1_j},...,p_{n_j},m_j,q_j,A_j) \geq -j^{-1}\,.
 \end{align*}
 Therefore we can find a maximizing sequence $(u_j,a_j,e_j) \in U\times\R^{d(n-i)}\times\R^d$ s.t.
  \begin{align*}
\left(-\partial_t\varphi+\big[\mathrm L^{u_j,a_j,e_j} +\mathrm F^{a_j,e_j} \big]\right)(t_j,z_j,p_{i+1_j},...,p_{n_j},m_j,q_j,A_j) \geq -2j^{-1}\,.
 \end{align*}
We define $b_j :=\frac1{\sqrt{1 + |a_{j}|^2+|e_j|^2}} \left (1, a_{j,i+1}^\top , \ldots , a_{j,n}^\top,  e_j^\top\right)^\top\in\mathcal{S}_{d(n-i)+d} \setminus \mathcal{D}_{d(n-i)+d}$, and get
 \begin{align*}
(b_{j_1})^2\left\{\left(-\partial_t\varphi+\big[\mathrm L^{u_j,\pr{\bar b}_j,\lf{\bar b_j}}+ \mathrm F^{\pr{\bar b}_j,\lf{\bar b_j}}\big]\right)(t_j,z_j,p_{i+1_j},...,p_{n_j},m_j,q_j,A_j) \right\}\geq -2j^{-1}(b_{j_1})^2\,.
\end{align*}
Appealing now to the relative compactness of the set $\mathcal{S}_{d(n-i)+d} \setminus \mathcal{D}_{d(n-i)+d}$ we obtain the existence
of a subsequence (still indexed by $j$) s.t. $\lim_{j\rightarrow\infty}b_j=\hat b$ with $\hat b\in\mathcal{S}_{d(n-i)+d}$. Moreover the compactness of $U$ ensures that, again up to a subsequence,  $\lim_{j\rightarrow\infty}u_j=\hat u\in U$. Therefore, using \eqref{eq condcv} and the continuity of the coefficients of $Z$, we obtain, after taking the limit over $j\to\infty$,
\begin{equation}\label{eq:super-H}
 \mathrm{H}^{\hat u,\hat b}\varphi( t, z,p_{i+1},..., p_n,m)\geq0\,,
 \end{equation}
 leading to the required result. By similar arguments we can prove the sub-solution property. In particular, for any smooth function  $\varphi$  s.t.
\begin{align}\label{eq:testsub}
(\text{strict})\max_{\interior(\mathscr{D}_i)}(w^*-\varphi)=(w^*-\varphi)(t,z, p_{i+1},..., p_n,m)=0\,,
\end{align}
one has
\begin{equation}\label{eq:sub-H}
\sup_{\substack{{(u,b)\in U\times\mathcal{S}_{d(n-i)+d}}}}\mathrm{H}^{u,b}\varphi( t, z, p_{i+1},..., p_n, m)\leq 0\,.
\end{equation}
The proof for  the choice of $\kappa$ and $\lambda$ identically equal to $1$ is then completed.\\
\noindent\textbf{Step 2.} We  extend the result to the case of general positive functions $\kappa$ and $\lambda$.
One can easily observe that \eqref{eq:super-H} (resp. \eqref{eq:sub-H}) is equivalent to writing that at point $(t, z,p_{i+1},..., p_n,m)$ and  for a smooth  function $\varphi$ satisfying \eqref{eq:testsuper} (resp. \eqref{eq:testsub}) one has
\begin{align}
 b^\top G(t,z, u,\partial_t \varphi,{\rm D} \varphi,{\rm D}^2 \varphi)  b\ge 0 , \quad \text{for some $u \in U,  b \in \mathcal{S}_{d(n-i)+d} $}\label{eq:Gsuper} \\
 \Big(\text{resp.}\quad b^\top G(t,z, u,\partial_t \varphi,{\rm D} \varphi,{\rm D}^2 \varphi) b \leq 0, \quad \text{for any $ u \in U, b \in \mathcal{S}_{d(n-i)+d} $} \Big),\label{eq:Gsub}
\end{align}
where for any $(t,z,u,c,q,A)\in [t_i,t_{i+1})\times \R^d\times U\times \R\times \R^{d+n-i+1}\times \mathbb S^{d+n-i+1}$ we denote by  $G(t,z,u,c,q,A)$  the following matrix in $\s^{1+d (n-i)+d}$,
\begin{gather*}
G(t,z,u,c,q,A)\defeq\left(\begin{matrix}  D & -\frac{1}{2}\left[\sigma^\top A^{zp_{i+1}} \right]^\top &\cdots&\cdots&\cdots&-\frac{1}{2}\left[\sigma^\top A^{zp_{n}} \right]^\top & -\frac{1}{2}\left[\sigma^\top A^{zm} \right]^\top \\
 -\frac{1}{2}\sigma^\top A^{p_{i+1}z^\top}& -\frac{1}{2}A^{p_{i+1},p_{i+1}}I_d&\textbf{0}&\cdots&\cdots&\textbf{0}& -\frac{1}{2}A^{p_{i+1}m}I_d \\
\vdots&\textbf{0}&\ddots&\ddots&\ddots&\vdots&\vdots\\
\vdots&\vdots&\ddots&\ddots&\ddots&\vdots&\vdots\\
 \vdots&\vdots&\ddots&\ddots&\ddots&\textbf{0}&\vdots\\
 -\frac{1}{2}\sigma^\top A^{p_{n}z^\top}&\textbf{0}&\cdots&\cdots&\textbf{0}&-\frac{1}{2}A^{p_{n},p_{n}}I_d  & -\frac{1}{2}A^{p_nm}I_d \\
 -\frac{1}{2}\sigma^\top A^{mz^\top}& -\frac{1}{2}A^{p_{i+1}m}I_d &\cdots&\cdots&\cdots& -\frac{1}{2}A^{p_{n}m}I_d&  -\frac{1}{2}A^{mm} I_d
\end{matrix}\right)\,,
\end{gather*}
with $\sigma\equiv\sigma(t,z,u)$ and $D\defeq-c-\mu^\top(t,z,u)q^z -\frac{1}{2}\Tr[\sigma\sigma^\top(t,z,u)A^{zz}]$.
Define the diagonal matrix
\begin{align*}
Q_{\kappa,\lambda}(p_{i+1},...,p_n,m)\defeq \text{diag}\left(1, \kappa(p_{i+1})I_d,...,\kappa(p_{n})I_d,\lambda(m)I_d\right)\,.
\end{align*}
A straightforward calculation shows that
\begin{align*}
& \det[(Q^{\top}_{\kappa,\lambda} GQ_{\kappa,\lambda})^{(k)}] = \lambda(m)^{2\max(k-(n-i)d-1,0)} \prod^{\lfloor\frac{k-1}{d}\rfloor +1}_{j=1} \kappa(p_{i+j})^{2\min(d,(1-j)d+k-1)}  \det[G^{(k)}]\,,
\end{align*}
where for any given matrix  $M\in \s^{1+d(n-i)+d} $ and $1\leq k\leq 1+d (n-i)+d$  we denote by $M^{(k)}\in\s^{k}$ the $k$-th leading principal sub-matrix of  $M$. Then, the positivity of the functions $\kappa$ and $\lambda$  implies that the quadratic form associated to $G$ and $Q^{\top}_{\kappa,\lambda} GQ_{\kappa,\lambda}$, respectively, have the same sign. Hence, it  follows from \eqref{eq:Gsuper} (resp.  \eqref{eq:Gsub}) that
\begin{align*}
 b^\top Q^{\top}_{\kappa,\lambda} G(t,z, u,\partial_t \varphi,{\rm D} \varphi,{\rm D}^2 \varphi) Q_{\kappa,\lambda}  b\ge 0 , \quad \text{for some $u \in U, b \in \mathcal{S}_{d(n-i)+d} $}\,,\\
 \Big(\text{resp.}\quad b^\top Q^{\top}_{\kappa,\lambda} G(t,z, u,\partial_t \varphi,{\rm D} \varphi,{\rm D}^2 \varphi) Q_{\kappa,\lambda}  b \leq 0, \quad \text{for any $ u \in U, b \in \mathcal{S}_{d(n-i)+d} $}\,, \Big)
\end{align*}
leading to
\begin{align*}
\sup_{\substack{{(u,b)\in U\times\mathcal{S}_{d(n-i)+d}}}}\mathrm{H}^{u,b}_{\kappa,\lambda}\varphi( t, z, p_{i+1},..., p_n, m)\geq 0\,\\
\Big( \text{resp.} \sup_{\substack{{(u,b)\in U\times\mathcal{S}_{d(n-i)+d}}}}\mathrm{H}^{u,b}_{\kappa,\lambda}\varphi( t, z, p_{i+1},..., p_n, m)\leq 0\Big)\,,
\end{align*}
which concludes the proof.
\qed

\subsection{On the space boundaries}\label{se bdary}
We study here the boundary conditions in $m$ and $p_k,\,1\le k\le n$.
We first divide $\partial\mathscr{D}_i,\,0\le i\le n-1,$ into different regions corresponding to the different boundaries associated with the level of controlled loss. More precisely, given $0\le i\le n-1$, we define
$\mathcal{P}_{i}\defeq\{I: I\subseteq \{i+1,...,n\},\,I\neq\emptyset\}.$ For any $I\in\mathcal{P}_{i}$, we  define $I^c\defeq\{i+1,\ldots, n\}\setminus I$, denote by $\text{Card}(I^c)$ its cardinality, and introduce $B_{i,I}\defeq\{(p_{i+1},...,p_n)\in[0,\infty)^{n-i}:\,p_{k}= 0\mbox{ for }k\in I\mbox{ and }p_k>0\mbox{ for }k\in I^c\},$ as well as $\mathscr{B}_{i,I}\defeq[t_i,t_{i+1})\times\R^{d}\times B_{i,I}.$ In particular,
$\mathscr{B}_{i}=\cup_{I\in\mathcal{P}_{i}}\mathscr{B}_{i,I}\cup \interior(\mathscr{B}_i).$
Then $\mathscr{D}_{i}=\interior(\mathscr{D}_i)\cup\partial \mathscr{D}_{i}$ where $\partial \mathscr{D}_{i}\defeq \left(\cup_{I\in\mathcal{P}_{i}}\mathscr{B}_{i,I} \times \R^+\right)\cup \left(\text{int}(\mathscr{B}_{i}) \times \{0\}\right).$
For any $0\leq i \leq n-1$, $I\in \mathcal{P}_{i}$, we define the following functions
\begin{itemize} \label{pg bdary}
\item on
 $\mathscr{C}_i$,
\begin{align*}
{w}_0(t,z,p_{i+1},...,p_n)&\defeq\inf_{(\nu,\alpha)\in \mathcal U\times \mathcal A^{n-i}}\E\left[f(Z^{t,z,\nu}_T)
+\sum^{n}_{k=i+1}\left(\Psi(Z^{t,z,\nu}_{t_k})-P^{t,p_{k},\alpha_k}_{t_k}\right)^+\right]\,.
\end{align*}
\item if $I^c\neq \emptyset$, on $[t_i,t_{i+1})\times\R^{d}\times\R^{\text{Card}(I^c)}\times\R$,
\begin{align*}
w_{1,I}(t,z,(p_k)_{k\in I^c},m)
&\defeq\inf_{\substack{\nu\in\mathcal U\\\alpha\in\mathcal A^{\text{Card}(I^c)}\\\eta\in\mathcal A}}
\E\left[\begin{matrix}\left(f(Z^{t,z,\nu}_T)-{M}^{t,m,\eta}_T\right)^+ +\sum_{k\in I}\Psi(Z^{t,z,\nu}_{t_k})
\\+\sum_{k\in I^c}\left(\Psi(Z^{t,z,\nu}_{t_k})-P^{t,p_{k},\alpha_{k}}_{t_k}\right)^+
\end{matrix}\right]\,,
\end{align*}
\item if $I^c=\emptyset$, on $[t_i,t_{i+1})\times\R^{d}\times\R$,
\begin{align*}
w_{1,I}(t,z,m)&\defeq\inf_{(\nu,\eta)\in\mathcal U\times \mathcal A}
\E\left[\left(f(Z^{t,z,\nu}_T)-{M}^{t,m,\eta}_T\right)^+ +\sum_{k=i+1}^n\Psi(Z^{t,z,\nu}_{t_k})\right]\,.
\end{align*}
\end{itemize}
We extend the definitions above to $t=T$ by setting $w_0(T,z) = f(z)$ on $\R^d$,  and $w_{1,I}(T,z,m)= (f(z)-m)^+$ on $\R^d\times \R$ for all $I\in\mathcal P_i$.
\medskip
\begin{remark}\label{re w1}
The functions $w_0$ and $(w_{1,I})_{I\in \mathcal P_i},\,0\le i\le n-1$ can be fully characterized on respectively, $\mathscr{B}_i$, $[t_i,t_{i+1})\times\R^{d}\times[0,\infty)^{\text{Card}(I^c)}\times\R^+$ if $I^c\neq \emptyset$ and $[t_i,t_{i+1})\times\R^{d}\times\R^+$ if $I^c= \emptyset$. This involves the same techniques as those developed to study the function $w$ (see above and hereinafter) but applied on a lower dimensional state space. For this reason, we consider from now on that these functions are a-priori known  and continuous on their respective domain of definition.
\end{remark}
The following proposition gives the natural Dirichlet conditions satisfied at the boundary $m=0$ and $p_k=0$, $i+1\le k\le n$.
\begin{proposition}\label{pr mbdary}
Fix $0\leq i\leq n-1$. On $\interior(\mathscr{B}_i)\times\{0\}$,
\begin{align}\label{eq BM}
w^*(t,z,p_{i+1},...,p_n,m)=w_*(t,z,p_{i+1},...,p_n,m)= w_0(t,z,p_{i+1},...,p_n)\,.
\end{align}
Moreover, on $\mathscr{B}_{i,I}\times\R^+,\,I\in \mathcal P_i$,
\begin{align}\label{eq BP}
w^*(t,z,p_{i+1},...,p_n,m)=w_*(t,z,p_{i+1},...,p_n,m)= w_{1,I}(t,z,(p_k)_{k\in I^c},m)\,.
\end{align}
\end{proposition}
\noindent{\it Proof}
We only prove \eqref{eq BM} as \eqref{eq BP} can be proved similarly. Fix $0\le i\le n-1$.
On one hand, by  the martingale property of stochastic integrals, one has on $\mathscr{C}_i\times\R$,
$$
w(t,z,p_{i+1},...,p_n,m)
\ge  w_0(t,z,p_{i+1},...,p_n)-m\,.
$$
On the other hand, on $\mathscr{C}_i\times\R$,
 \begin{align*}
& w(t,z,p_{i+1},...,p_n,m)
\le\inf_{(\nu,\alpha)\in \mathcal{U}\times \mathcal{A}^{n-i}}\E\left[\left(f(Z^{t,z,\nu}_T)-M^{t,m,0}_T\right)^+ +\sum^{n}_{k=i+1}\left(\Psi(Z^{t,z,\nu}_{t_k})-P^{t,p_{k},\alpha_{k}}_{t_k}\right)^+\right]\,.
\end{align*}
As $f\geq 0$, the latter leads to $w(t,z,p_{i+1},...,p_n,m)\le w_0(t,z,p_{i+1},...,p_n)$
 when $m\ge0$. We conclude by taking the upper and lower limit and recalling Remark \ref{re w1}.
\qed
\subsection{A comparison principle for \eqref{viscsolBETbis}}\label{se CP}
We consider ${\lambda}:m\in\R\mapsto \lambda(m)\defeq 1\vee m>0$ and ${\kappa}:p\in\R\mapsto \kappa(p)\defeq 1\vee p>0$. The operator in \eqref{viscsolBETbis} is non-standard as it involves a non-linearity in the time-derivative. However, thanks to a strict super-solution approach (see e.g. \cite{IL90,CST05}), a comparison result can be proved.
\begin{lemma}[Strict super-solution property]\label{le StrictSupersol}
Fix $0\le i\le n-1$. Let us define on $[t_i,t_{i+1})\times\R^{n-i}\times\R$ the smooth positive function
$$\phi(t,p_{i+1},...,p_n,m) \defeq e^{(t_{i+1}-t)} \left(1+\sum^n_{k=i+1}\ln(1+p_k)+\ln(1+m)\right)\,.$$
Let $v$ be a lower semi-continuous viscosity super-solution of \eqref{viscsolBETbis}.
Then the function $v+\xi \phi$, $\xi>0$, satisfies in the viscosity sense
\begin{align}
 \sup_{\substack{{(u,b)\in U\times\mathcal{S}_{d(n-i)+d}}}} \mathrm H_{\kappa,\lambda}^{u,b}(v+\xi \phi) \geq \xi \frac{1}{8}\text{ on }\interior(\mathscr{D}_i)\,. \label{eq strict supersol}
\end{align}
\end{lemma}
\noindent{\it Proof}
Fix $0\le i\le n-1$ and $(t, z, p_{i+1},..., p_n, m)\in \interior(\mathscr{D}_i)$.  Let $\xi>0$ and $\varphi$ be a smooth function s.t. $\min_{\interior(\mathscr{D}_i)}((v+\xi \phi)-\varphi)=((v+\xi \phi)-\varphi)( t, z, p_{i+1},..., p_n, m)=0$.
Since $ \phi$ is a smooth function, the function $\psi\defeq\varphi-\xi \phi$ is a test function for $v$.

Let $b \in \mathcal{S}_{d(n-i)+d} \setminus \mathcal{D}_{d(n-i)+d}$ and $u\in U$. We obtain by definition of $\mathrm H_{\kappa,\lambda}^{u,b}$,
\begin{align}\label{eq ineqH}
\mathrm H_{\kappa,\lambda}^{u,b}\varphi( t, z, p_{i+1},..., p_n, m)\ge \mathrm H_{\kappa,\lambda}^{u,b}\psi( t, z, p_{i+1},..., p_n, m) + \mathfrak{A} \,,
\end{align}
where
\begin{align*}
\mathfrak{A}= \xi (b_1)^2\left(\begin{matrix} -\partial_t \phi(\cdot) - \frac12 \sum^n_{k=i+1} \kappa( p_k)^2|\pr{\bar b}_k|^2 \mathrm D_{p_k p_k} \phi(\cdot) - \frac12 \lambda( m)^2|\lf{\bar b}|^2 \mathrm D_{mm} \phi(\cdot)\end{matrix}\right)( t, z, p_{i+1},..., p_n, m)\,.
\end{align*}
We now provide a lower bound for $\mathfrak{A}$ and compute
\begin{align*}
 \mathfrak{A}=  &\xi (b_1)^2  e^{(t_{i+1}- t)} \left(\begin{matrix} 1+ \sum^n_{k=i+1}\ln(1+ p_k)+\ln(1+ m) \end{matrix}\right)\\&+ \xi (b_1)^2  e^{(t_{i+1}- t)} \left(\begin{matrix} \sum^n_{k=i+1}\frac{\kappa(p_k)^2|\pr{\bar b}_k|^2}{2(1+ p_k)^2}+\frac{\lambda(m)^2|\lf{\bar b}|^2}{2(1+ m)^2} \end{matrix}\right)
\ge  \frac{\xi(b_1)^2 }{8} \left( 1+ |\pr{\bar b}|^2+|\lf{\bar b}|^2\right)\,,
\end{align*}
since $ m,p_k\ge 0,\,i+1\le k\le n,$ and for any $l\ge 0$, $\frac{1\vee l^2}{(1+l)^2}\ge \frac{1}{4}.$
Noticing that $(b_1)^2(1+|\pr{\bar b}|^2+|\lf{\bar b}|^2)=1$, we finally obtain
$
\mathfrak{A} \ge\xi\frac{1}{8}.
$
Thanks to the arbitrariness of $u$ and $b$ and \eqref{eq ineqH}, one has, after appealing to the super-solution property of $v$,
\begin{align*}
 \sup_{\substack{{(u,b)\in U\times\mathcal{S}_{d(n-i)+d}\setminus \mathcal{D}_{d(n-i)+d}}}} \mathrm H_{\kappa,\lambda}^{u,b} \varphi(t, z, p_{i+1},..., p_n, m)\ge \xi \frac{1}{8}\,.
\end{align*}
We finally conclude the proof recalling Remark \ref{re opcont}.
\qed

We can now state a comparison result holding for viscosity solutions of \eqref{viscsolBETbis} whose proof is postponed to Appendix B and is based on Lemma \ref{le StrictSupersol}.
\begin{theorem}[Comparison principle]\label{CparisonPrinciple}
Fix $0\le i\le n-1$. Let $V$ (resp. $U$) be a lower semi-continuous (resp. upper semi-continuous) function  satisfying
$$
\left| V(t,z,p_{i+1},...,p_n,m)\right| + \left| U(t,z,p_{i+1},...,p_n,m)\right| \leq C\left(1+|z|\right)\text{ on }\mathscr{D}_i\,.
$$
Moreover, assume that
on $\interior(\mathscr{D}_i)$, $V$ (resp. $U$) is a viscosity super-solution (resp. sub-solution) of \eqref{viscsolBETbis},
 on $\partial \mathscr{D}_i$, $V(\cdot)\ge U(\cdot)$, and
 on $\R^{d}\times[0,\infty)^{n-i}\times\R^+$, $V(t_{i+1},\cdot) \ge U(t_{i+1},\cdot)$.
 Then $V\geq U$ on $\mathscr{D}_i$.
\end{theorem}
\subsection{A complete characterization of $w$}\label{se cpletecharact}
Thanks to the results in the previous sections we can now obtain a full characterization of $w$ by the HJB equation. Moreover, we obtain the time-continuity of $w$ on each interval, which completes the result derived in Proposition \ref{prop growth}. Indeed, unlike \cite{BPZ16}, uniform $L^2$-boundedness of the admissible controls is not considered and the regularity of $w$ cannot be proven a priori.
\begin{theorem}[Complete characterization of $w$]\label{th charactw}
The function $w$  is the unique viscosity solution of \eqref{viscsolBETbis} on $\interior(\mathscr{D}_i)$  for any $0\le i\le n-1$, in the class of functions being continuous on $\mathscr D_i$ and satisfying the growth condition \eqref{eq polygrowth} and with the following terminal and boundary conditions
$$
{w}(T,z,m)=(f(z)-m)^+\; \mbox{ on } \R^{d}\times\R^+,
$$
\begin{align}\label{eq:mp}
& w={w}_0\;\mbox{ on } \interior(\mathscr{B}_i)\times\{0\},\qquad {w}={w}_{1,I}\; \mbox{ on }\mathscr{B}_{i,I}\times\R^+,\,\forall I\in \mathcal P_i\,,
\end{align}
 and
\begin{align}\label{eq:tip1}
 \lim_{t\uparrow t_{i+1}}{w}(t,z,p_{i+1},...,p_n,m) = {w}(t_{i+1},z,p_{i+2},...,p_n,m)+ (\Psi(z)-p_{i+1})^+ \; \mbox{ on } \R^{d}\times[0,\infty)^{n-i}\times\R^+.
 \end{align}
\end{theorem}
\noindent{\it Proof}
By definition, the condition on $\{T\}\times\R^d\times\R$ is satisfied. Additionally, we know from Proposition \ref{prop growth} (resp. Proposition \ref{pr mbdary} and Remark \ref{re w1} ) that $w$  satisfies the linear growth condition \eqref{eq polygrowth} (resp. the boundary conditions in \eqref{eq:mp} and is continuous on $\left(\interior(\mathscr{B}_i)\times\{0\}\cup (\cup_{I\in \mathcal P_i}\mathscr{B}_{i,I}\times\R^+)\right)$ for any $0\le i\le n-1$). Moreover, it follows from Theorem  \ref{theo caractinside} that for any $0\leq i\leq n-1$,  $w^*$ (resp. $w_*$) is an upper semi-continuous  (resp. lower semi-continuous) viscosity sub-solution (resp. super-solution) to \eqref{viscsolBETbis} on $\interior(\mathscr{D}_i)$. To prove the continuity property on $\mathscr{D}_i,\,1\le i\le n-1$, we proceed by induction on $i$. The uniqueness property and \eqref{eq:tip1} are a by-product of the proof by induction. \\
Let $i=n-1$. Thanks to Proposition \ref{prop growth} and \ref{pr mbdary} as well as Theorem  \ref{theo caractinside}, the uniqueness of the solution to \eqref{viscsolBETbis} and continuity of $w$ on $\interior(\mathscr{D}_{n-1})$ follow from Theorem \ref{CparisonPrinciple}.
Let us now assume that $w$ is continuous on $\interior(\mathscr D_{i+1})$ for some $0\le i\le n-2$, and show its continuity on $\interior(\mathscr D_{i})$.
The result follows by the same arguments as above once proved \eqref{eq:tip1}, in virtue of the Lipschitz continuity of $w$ in the space variables (remind Proposition \ref{prop growth}).
To this aim we start by introducing on $[0,t_{i+2})\times\R^d\times\R^{n-i-1}\times\R,\,0\le i\le n-2$, the auxiliary function
\begin{align*}
&\hat w(t,z,p_{i+2},...,p_n,m)=\inf_{\substack{{(\nu,\alpha,\eta)\in\mathcal{U}\times {\mathcal{A}}^{n-i-1}\times\mathcal{A}}}}{\mathrm{J}}^{{\nu},\alpha,\eta}(t,z,p_{i+2},...,p_n,m)\,.
\end{align*}
We observe that on $[t_{i+1},t_{i+2})\times\R^d\times\R^{n-i-1}\times\R,\,0\le i\le n-2$,
\begin{align*}
&\hat w(t,z,p_{i+2},...,p_n,m)= w(t,z,p_{i+2},...,p_n,m)\,.
\end{align*}
Moreover, $\hat w^*$ and $\hat w_*$ satisfy a linear growth condition on their respective domain. Additionally, the induction assumption implies that $\hat w$ is continuous on $\mathscr D_{i+1}$. Hence
$$\hat w^*(t_{i+2},z,p_{i+2},...,p_n,m)=\hat w_*(t_{i+2},z,p_{i+2},...,p_n,m)\mbox{ on }\R^d\times\R^{n-i-1}\times\R\,.$$ Moreover, proceeding as in Propositions \ref{pr mbdary} (resp. Theorem  \ref{theo caractinside}), one can derive continuous boundary conditions on $[0,t_{i+2})\times\R^d\times(0,\infty)^{n-i-1}\times\{0\}$ and $[0,t_{i+2})\times\R^d\times B_{i+1,I}\times\R^+\,,\forall I\in \mathcal P_{i+1}$ (resp. characterize $\hat w^*$ and $\hat w_*$ on $[0,t_{i+2})\times\R^d\times(0,\infty)^{n-i-1}\times\R^+_*$). Therefore, appealing to Theorem \ref{CparisonPrinciple}, one obtains the continuity of $\hat w$ on $[0,t_{i+2})\times\R^d\times(0,\infty)^{n-i-1}\times\R^+_*$.

As a result, for any $(z,p_{i+2},\ldots,p_n,m)\in \R^d\times [0,+\infty)^{n-i-1}\times \R^+$, one has
\begin{align}\label{eq:ip1}
 \lim_{h\to 0} \left|\hat w(t_{i+1}-h,z,p_{i+2},\ldots,p_n,m) - w(t_{i+1},z,p_{i+2},\ldots,p_n,m)\right| = 0\,.
\end{align}
Let $h>0$ be  s.t. $t_{i+1}-h\in [t_{i}, t_{i+1})$. On $\R^d\times \R^{n-i}\times \R$, one can easily check that
\begin{gather*}
 w(t_{i+1}-h,z,p_{i+1},\ldots, p_n,m)-\hat w(t_{i+1}-h,z,p_{i+2},\ldots, p_n,m)
\leq \sup_{\nu\in \mathcal U}\mathbb E\left[\Big(\Psi(Z^{t_{i+1}-h, z,\nu}_{t_{i+1}}) - P^{t_{i+1}-h,p_{i+1}, 0}_{t_{i+1}}\Big)^+\right]\,.
\end{gather*}
Moreover, from the martingale property of stochastic integrals follows
\begin{gather*}
 {w}(t_{i+1}-h,z,p_{i+1},\ldots, p_n,m)  -\hat w(t_{i+1}-h,z,p_{i+2},\ldots, p_n,m)
  \geq  \Big(\inf_{\nu\in \mathcal U}\mathbb E\left[\Psi( Z^{t_{i+1}-h,z,\nu}_{t_{i+1}})\right] -p_{i+1} \Big)^+ \,.
\end{gather*}
Therefore, thanks to the Lipschitz continuity of  $\Psi$ together with  classical estimates on the process $Z$, there exists a uniform $C>0$ s.t.
\begin{gather*}
|w(t_{i+1}-h,z,p_{i+2},\ldots, p_n,m) -\hat w(t_{i+1}-h,z,p_{i+2},\ldots, p_n,m)- (\Psi(z)-p_{i+1})^+|\le C \sqrt h \left(1+|z|\right)\,.
\end{gather*}
Condition \eqref{eq:tip1} then follows by taking the limit when $h$ tends to zero and after appealing to \eqref{eq:ip1}.
\qed
\begin{remark}
The techniques developed in this paper can be applied to the case of \textit{next-period} controlled-loss constraints, with the terminology of \cite{JKT17}. This type of constraints has been studied for stochastic target problems under a complete market setting in \cite{B18b}. For a given $\nu\in\mathcal{U}$, they write on $\mathscr{C}_i,\, 0\le i\le n-1$, as $\E\left[\Psi( Z^{t,z,\nu}_{t_{i+1}})|\mathcal{F}_{t}\right]\le p_{i+1}$ and $\E\left[\Psi( Z^{t,z,\nu}_{t_k})|\mathcal{F}_{t_{k-1}}\right]\le p_k,\,i+2\le k\le n$, and appealing to the techniques developed in Section \ref{se PbRdPbedp}, one can prove that they re-write as
$\Psi(Z^{t,z,\nu}_{t_{i+1}} )\le P^{t,p_{i+1},\alpha_{i+1}}_{t_{i+1}}$ and $\Psi(Z^{t,z,\nu}_{t_k} )\le P^{t_{k-1},p_{k},\alpha_{k}}_{t_{k}},\,i+2\le k\le n,$ for some $(\alpha_{i+1},...,\alpha_n,\eta)\in\mathcal{A}^{n-i}\times\mathcal{A}$.
\end{remark}
\section{Conclusions}
Assuming the existence of an optimizer for the value function $w$ defined in \eqref{defw}, we  proved that the original value function $V$ in \eqref{pbintialedp} can be described as the zero level of $w$. This result has the great advantage of providing a characterization of the value function associated with a state constrained optimal control problem without requiring any viability or strong regularity assumptions on the coefficients of the diffusion process. However, $w$  is associated with an (unconstrained) optimal control problem involving unbounded controls, which raises additional difficulties in the treatment of the associated PDE.
Here, we  provide a full characterization of the level set function $w$ as the unique piecewise continuous viscosity solution of a suitable HJB equation passing through a compactification of the differential operator.\\
\indent This paper opens new avenues for further research which includes  the treatment of the case where the constraints hold in probability, and the study of the numerical approximation of $V$ appealing to the characterization of $w$.

\begin{acknowledgements}
Authors are grateful to the Young Investigator Training Program and Association of Bank Foundations. This is a pre-print of an article published in Journal of Optimization Theory and Applications. The final authenticated version is available online at: https://doi.org/10.1007/s10957-020-01724-8
\end{acknowledgements}
\appendix
\section*{Appendix A}\label{DPP}
\noindent\it{Proof of Theorem \ref{th dpp} }
Fix $(t,z,p_{i+1},...,p_n,m)\in\mathscr{C}_i\times\R,\,0\le i\le n-1$, and a stopping time $t\le \theta< t_{i+1}$. We denote $\widehat w$ the right-hand side of \eqref{eq dpp}. \\
\textbf{Step 1. Proof of $w\ge \widehat w$.} By the definition of $w$ in \eqref{defw} and the Flow property one has
\begin{gather*}
w(t,z,p_{i+1},...,p_n,m)\ge \inf_{(\nu,\alpha,\eta)\in\mathcal{U}\times\mathcal{A}^{n-i}\times\mathcal{A}}\E\left[ w\left({\theta},{Z}^{t,z,{\nu}}_{\theta},{P}^{t,p_{i+1},\alpha_{i+1}}_{\theta},...,{P}^{t,p_{n},\alpha_n}_{\theta},
{M}^{t,m,\eta}_{\theta}\right)\right]\,,
\end{gather*}
leading to $w(t,z,p_{i+1},...,p_n,m)\ge  \widehat w(t,z,p_{i+1},...,p_n,m).$\\
\textbf{Step 2. Proof of $w\le \widehat w$.}
  We fix $(\hat{{\nu}},\hat\alpha_{i+1},...,\hat\alpha_{n},\hat\eta)\in{\mathcal{U}}\times \mathcal{A}^{n-i}\times\mathcal{A}$, and consider $\mu$, the measure induced by $(\theta,\xi,\zeta_{i+1},...,\zeta_n,\kappa)$  on $\mathscr{C}_i\times\R$ with $(\xi,(\zeta_k)_{i+1\le k\le n},\kappa)\defeq ({Z}^{t,z,\hat{{\nu}}}_{{\theta}},({P}^{t,p_{k},\hat{\alpha}_k}_{{\theta}})_{i+1\le k\le n},{M}^{t,m,\hat{\eta}}_{{\theta}})$.  We appeal to \cite[Proposition 7.50, Lemma 7.27]{BS78} to prove that, for each $\varepsilon>0$, we can build $(n-i+2)$ Borel-measurable maps ${\nu}^\varepsilon_\mu$, $\alpha^\varepsilon_\mu\equiv({\alpha}_{i+1,\mu}^{\varepsilon},...,{\alpha}_{n,\mu}^{\varepsilon})$ and ${\eta}_\mu^{\varepsilon}$ s.t.
  $({\nu}^\varepsilon_\mu,\alpha^\varepsilon_\mu,{\eta}_\mu^{\varepsilon})\in{\mathcal{U}}\times\mathcal{A}^{n-i}\times\mathcal{A},$ and
   \begin{align}\label{J}
 &w(\theta,\xi,\zeta_{i+1},...,\zeta_n,\kappa)\ge {\mathrm{J}}^{{\nu}^\varepsilon_\mu,\alpha^\varepsilon_\mu,{\eta}_\mu^{\varepsilon}}(\theta,\xi,\zeta_{i+1},...,\zeta_n,\kappa)
-\varepsilon\,.
\end{align}
We now use \cite[Lemma 2.1]{ST02GEO} to obtain $\nu^{\varepsilon}$, $\alpha^\varepsilon$ and ${\eta}^{\varepsilon}$ s.t.
\begin{align*}
\nu^{\varepsilon}\1_{[\theta,T]}&=\nu^{\varepsilon}_\mu(\theta,\xi,\zeta_{i+1},...,\zeta_n,\kappa)\1_{[\theta,T]}\,\mathrm{d}t\times\mathrm{d}\p\text{-a.e.}\,,
\\{\alpha}_{k}^{\varepsilon}\1_{[\theta,t_{k}]}&={\alpha}_{k,\mu}^{\varepsilon}(\theta,\xi,\zeta_{i+1},...,\zeta_n,\kappa)\1_{[\theta,t_{k}]}\,
\mathrm{d}t\times\mathrm{d}\p\text{-a.e.}, i+1\le k\le n\,,
\\\eta^{\varepsilon}\1_{[\theta,T]}&=\eta^{\varepsilon}_\mu(\theta,\xi,\zeta_{i+1},...,\zeta_n,\kappa)\1_{[\theta,T]}\,\mathrm{d}t\times\mathrm{d}\p\text{-a.e.}
\end{align*}
This implies that  $\hat{\nu}^\varepsilon\defeq\hat\nu\1_{[t,\theta)}+\nu^{\varepsilon}\1_{[\theta,T]}\in{\mathcal{U}}$, $\hat{\alpha}_{k}^{\varepsilon}\defeq\hat{\alpha}_{k}\1_{[t,\theta)}+\alpha^{\varepsilon}_{k}\1_{[\theta,t_{k}]}\in \mathcal A,\, i+1\leq k\leq n$ and $\hat{\eta}^{\varepsilon}\defeq\hat{\eta}\1_{[t,\theta)}+\eta^{\varepsilon}\1_{[\theta,T]}\in \mathcal A$ and
 \eqref{J} holds where,
 according to \cite[Remark 6.1]{BEI10},
\begin{align*}
&{\mathrm{J}}^{({\nu}^\varepsilon_\mu,\alpha^\varepsilon_\mu,{\eta}_\mu^{\varepsilon})}(\theta,\xi,\zeta_{i+1},...,\zeta_n,\kappa)
=\E\left[\left.
\left(f({Z}^{\theta,\xi,\nu^{\varepsilon}}_T)-{M}^{\theta,\kappa,\eta^{\varepsilon}}_T\right)^+
+\sum_{k=i+1}^n\left(\Psi({Z}^{\theta,\xi,\nu^{\varepsilon}}_{t_{k}})-{P}^{\theta,\zeta_{k},\alpha^{\varepsilon}_{k}}_{t_{k}}\right)^+
\right|\left(\theta,\xi,\zeta_{i+1},...,\zeta_n,\kappa\right)\right]
\,.
\end{align*}
We conclude taking the expectation on both sides in \eqref{J}, using the tower property of expectation, and appealing to the arbitrariness of $(\hat{{\nu}},\hat\alpha_{i+1},...,\hat\alpha_{n},\hat\eta)\in{\mathcal{U}}\times \mathcal{A}^{n-i}\times\mathcal{A}$ and  $\varepsilon$.\qed
\section*{Appendix B}\label{app:CP}
We first state the following lemma which is involved in the proof of Theorem \ref{CparisonPrinciple}.
\begin{lemma}[Modulus of continuity]\label{le ModulusCont0} Fix $0\le i\le n-1$. There exists $\rho>0$ s.t.:
for any  $t\in [t_i,t_{i+1})$, $z,r\in \R^d$, $m,l\in \R^+$, $p,q\in [0,\infty)^{n-i}$ (with  $p\defeq(p_{i+1},...,p_n)^\top$ and $q\defeq(q_{i+1},...,q_n)^\top$), for any $\cX, \cY\in\s^{d+n-i+1}$ satisfying
\begin{align}\label{A10}
\left(\begin{matrix} \cX&0\\0&-\cY \end{matrix}\right)&\leq
\frac{3}{\varepsilon}\left(\begin{matrix} I &-I\\-I& I \end{matrix}\right)+2\zeta\text{e}^{-\rho t}\left(\begin{matrix}\bar I &\textbf{0}\\\textbf{0}& \bar I \end{matrix}\right)
\,,
\end{align}
for $\zeta,\varepsilon>0$, with $I\equiv I_{d+n-i+1}$, and $\bar I\defeq
\mathrm{diag} (I_{d}, 0, \ldots, 0)\in\s^{d+n-i+1}$,
 and for
\begin{align*}
& c_1,c_2\mbox{ s.t. } c_1-c_2= -\zeta\rho \text{e}^{-\rho t}\left(1+|z|^{2}\right)- \zeta\rho \text{e}^{-\rho t}\left(1+|r|^{2}\right)\in\R^-\,,\\
& \Delta_1\defeq \left(\begin{matrix}\frac{1}\varepsilon (z-r)+2\zeta\text{e}^{-\rho t}z\\\frac{1}\varepsilon (p-q)\\\frac{1}\varepsilon (m-l)  \end{matrix} \right)\,,\, \,\,\Delta_2\defeq \left(\begin{matrix}\frac{1}\varepsilon (z-r)-2\zeta\text{e}^{-\rho t} r\\\frac{1}\varepsilon (p-q)\\\frac{1}\varepsilon (m-l)\end{matrix} \right)\,\in \R^{d+n-i+1}\,,\\
\end{align*}
one has, with $\Theta_1\defeq(t,z,p,m,\Delta_1,\cX,c_1)$ and $\Theta_2\defeq(t,r,q,l,\Delta_2,\cY,c_2)$,
\begin{align*}
\sup_{\substack{{(u,b)\in U\times\mathcal{S}_{d(n-i)+d}}}} \hspace{-0.cm}\mathrm H_{\kappa,\lambda}^{u,b} (\Theta_2)- \sup_{\substack{{(u,b)\in U\times\mathcal{S}_{d(n-i)+d}}}} \hspace{-0.cm}\mathrm H_{\kappa,\lambda}^{u,b}(\Theta_1)
\leq\frac{C}{\varepsilon}\left(|z-r|^2+|p-q|^2+(m-l)^2\right)\,,
\end{align*}
for some $C>0$.
\end{lemma}
\noindent{\it Proof}
Consider $\Theta_1$ and $\Theta_2$ defined in the theorem. We notice (recall Remark \ref{re opcont}),
\begin{align*}
& \sup_{\substack{{(u,b)\in U\times\mathcal{S}_{d(n-i)+d}}}}\mathrm H_{\kappa,\lambda}^{u,b} (\Theta_2)-\sup_{\substack{{(u,b)\in U\times\mathcal{S}_{d(n-i)+d}}}}\mathrm H_{\kappa,\lambda}^{u,b}(\Theta_1)
\\&\quad\quad\quad\quad\quad\quad\quad\quad\quad\quad\quad\quad\quad\quad\le  \sup_{\substack{{u\in U}\\{b\in\cS_{d(n-i)+d}\setminus\cD_{d(n-i)+d}}}}\left\{ \mathrm H_{\kappa,\lambda}^{u,b} (\Theta_2)- \mathrm H_{\kappa,\lambda}^{u,b} (\Theta_1)\right\}\;.
\end{align*}
For $b \in \cS_{d(n-i)+d}\setminus\cD_{d(n-i)+d}$ and $u\in U$, we compute by definition of $\mathrm H_{\kappa,\lambda}^{u,b}$, \begin{align*}
\mathrm H_{\kappa,\lambda}^{u,b}(\Theta_2)-\mathrm H_{\kappa,\lambda}^{u,b}(\Theta_1) \le (b_1)^2\left(\mathfrak{A} + \mathfrak{B}+ \mathfrak{C} \right)\,,
\end{align*}
where
\begin{align*}
\mathfrak{A} &= \frac{1}{\varepsilon}\left(\mu(t,z,u)- \mu(t,r,u)\right)^\top\left(z-r\right),
\end{align*}
\begin{align*}
\mathfrak{B} &= 2\zeta\text{e}^{-\rho t}\mu^\top(t,z,u)z+2\zeta\text{e}^{-\rho t}\mu^\top(t,r,u)r-\zeta\rho\text{e}^{-\rho t}\left(2+|z|^{2}+|r|^{2}\right)\,,
\end{align*}
and
\begin{align*}
\mathfrak{C} = -\frac12 \Tr\left[\bar{\sigma}\bar{\sigma}^\top(t,r,q,l,u,b)\cY \right ] + \frac12 \Tr \left [\bar{\sigma}\bar{\sigma}^\top(t,z,p,m,u,b)\cX \right ]\,,
\end{align*}
where for  any $(t,r,q,l)\in[t_i,t_{i+1})\times\R^{d}\times[0,\infty)^{n-i}\times\R^+$, $b\in\cS_{d(n-i)+d}\setminus\cD_{d(n-i)+d}$ and $u\in U$,
$\bar{\sigma}(t,r,q,l,u,b)\defeq \left(\begin{matrix}\sigma(t,r,u),\kappa(q_{i+1})\bar b^{\flat\,\top}_{i+1},\cdots,\kappa(q_{n})\bar b^{\flat\,\top}_n,\lambda(l)\bar b^{\sharp\,\top}\end{matrix} \right)^\top.$

Using the Lipschitz and growth properties of $\mu$, we obtain some $C, \hat C>0$ s.t.
\begin{align*}
\mathfrak{A} \le \frac{ C}{\varepsilon}|z-r|^2\quad
\mbox{and}\quad
\mathfrak{B} \le \zeta \hat C\text{e}^{-\rho t}\left(1+|z|^{2}+|r|^{2}\right)
-\zeta\rho\text{e}^{-\rho t}\left(1+|z|^{2}+|r|^{2}\right)
\,.
\end{align*}
For $\mathfrak{C}$, we use \eqref{A10} and the Lipschitz continuity of $\sigma$, $\kappa$ and $\lambda$ to get some $\bar C>0$ s.t.
\begin{align*}
\mathfrak{C} \le& \bar C\left(\begin{matrix}\frac{1}{\varepsilon}\left(1+|\pr{\bar b}|^2+|\lf{\bar b}|^2\right)\left(|z-r|^2+|p-q|^2+(m-l)^2\right)
+\zeta\text{e}^{-\rho t}\left(1+|z|^{2}+|r|^{2}\right)\end{matrix}\right)
\,.
\end{align*}
Taking $\rho\ge \hat C+\bar C+1$ for instance, we obtain for some $C>0$,
$$\mathfrak{B}+\mathfrak{C}\le  \frac{C}{\varepsilon}\left(1+|\pr{\bar b}|^2+|\lf{\bar b}|^2\right)\left(|z-r|^2+|p-q|^2+(m-l)^2\right)
\,.$$
The proof is concluded by observing that $(b_1)^2(1+|\pr{\bar b}|^2+|\lf{\bar b}|^2)= 1$.
\qed

We can now prove Theorem \ref{CparisonPrinciple}.

\noindent\it{Proof of Theorem \ref{CparisonPrinciple}.}
Fix $0\le i\le n-1$. For $\xi>0$, we introduce on $(t_i,t_{i+1})\times\R^{d}\times[0,\infty)^{n-i}\times\R^+$ the following auxiliary functions
\begin{align*}
V_{\xi}(t,z,p,m)\defeq (V+\xi \phi)(t,z,p,m) + \xi\left(\frac{1}{t-t_i}\right),\quad
U_{\xi}(t,z,p,m)\defeq (U-\xi \phi)(t,z,p,m)\,,
\end{align*}
where $p\defeq(p_{i+1},...,p_n)^\top$ and
with $\phi$ defined in Lemma \ref{le StrictSupersol}. Appealing to Lemma \ref{le StrictSupersol}, one can easily check that $V_{\xi}$ is a strict super-solution
 of \eqref{viscsolBETbis} satisfying \eqref{eq strict supersol}. Analogously $U_\xi$ can be proved to be a sub-solution of \eqref{viscsolBETbis}.

We prove that $U-V\leq0$ on $\mathscr{D}_i$. To this aim we first show arguing by contradiction that for all $\xi>0$, $(U_{\xi}-V_{\xi})\leq0$ on $(t_i,t_{i+1})\times\R^{d}\times[0,\infty)^{n-i}\times\R^+$, and the proof is completed sending $\xi$ to zero.\\
\textbf{Step 1.}
We assume to the contrary that we can find $\xi>0$ s.t.
\begin{align}\label{eq CPinit}
\sup_{(t_i,t_{i+1})\times\R^{d}\times[0,\infty)^{n-i}\times\R^+} (U_{\xi}-V_{\xi})>0\,.
\end{align}
We define on $(t_i,t_{i+1})\times\R^{d}\times[0,\infty)^{n-i}\times\R^+$,
$$\Phi_{\xi,\zeta}(t,z,p,m)\defeq(U_\xi - V_\xi)(t,z,p,m) - 2\zeta e^{-\rho t}  (1+|z|^{2})\,,$$
for $\zeta>0$ and with $\rho>0$ defined in Lemma \ref{le ModulusCont0}.
Using the growth conditions and semi-continuity of $U$ and $V$ as well as \eqref{eq CPinit} we obtain that for $\xi,\zeta>0$ small enough
\begin{align*}
0<M\defeq\sup_{(t_i,t_{i+1})\times\R^{d}\times[0,\infty)^{n-i}\times\R^+} \Phi_{\xi,\zeta}(t,z,p,m)<\infty\,.
\end{align*}
On $(t_i,t_{i+1})\times\R^{d}\times\R^{d}\times[0,\infty)^{(n-i)}\times[0,\infty)^{(n-i)}\times\R^+\times\R^+,$ set
\begin{align*}
\Psi_{\xi,\zeta,\varepsilon}(t,z,r,p,q,m,l) \defeq & U_{\xi}(t,z,p,m)-V_{\xi}(t,r,q,l)
-\zeta\text{e}^{-\rho t}(1+|z|^{2})-\zeta\text{e}^{-\rho t}(1+|r|^{2})\\
&-\frac{1}{2\varepsilon}\left(|z-r|^{2} +|p-q|^2+(m-l)^2\right)\,,
\end{align*}
for $\varepsilon>0$ and with $q\defeq(q_{i+1},...,q_n)^\top$.
Again, the growth conditions and semi-continuity of  $U$ and $V$ ensure that for $\xi,\zeta,\varepsilon>0$ the function $\Psi_{\xi,\zeta,\varepsilon}$ admits a maximum $M_\varepsilon$ at $(t_{\varepsilon},z_{\varepsilon},r_{\varepsilon},p_{\varepsilon},q_{\varepsilon},m_{\varepsilon},l_{\varepsilon})$, with $p_\varepsilon\defeq( p_{i+1_\varepsilon},..., p_{n_\varepsilon})^\top$ and $q_\varepsilon\defeq( q_{i+1_\varepsilon},..., q_{n_\varepsilon})^\top$, on $(t_i,t_{i+1})\times\R^{d}\times\R^{d}\times[0,\infty)^{(n-i)}\times[0,\infty)^{(n-i)}\times\R^+\times\R^+$ (we omit the dependency on $(\xi,\zeta)$ for the sake of clarity).
Using standard arguments (see, for instance, the proof of \cite[Theorem 4.4.4]{P09} and \cite[Lemma 3.1]{CIL92}), one can prove that there exists $(\bar t,\bar z,\bar p,\bar m)\in (t_i,t_{i+1})\times \R^{d}\times [0,\infty)^{n-i}\times \R^+$, with $\bar p\defeq(\bar p_{i+1},...,\bar p_n)^\top$, s.t.
\begin{align}
\begin{cases}\label{eq cv}
\lim_{\varepsilon \downarrow 0} t_\varepsilon= \bar t,\,\lim_{\varepsilon \downarrow 0} z_\varepsilon,r_\varepsilon= \bar z,\,\lim_{\varepsilon \downarrow 0}  p_\varepsilon,q_\varepsilon= \bar p,\,\lim_{\varepsilon \downarrow 0}  m_\varepsilon,l_\varepsilon= \bar m\,,
\\
\lim_{\varepsilon \downarrow 0} \frac{1}{\varepsilon}\left(|z_\varepsilon-r_\varepsilon|^{2} +|p_\varepsilon-q_\varepsilon|^2+(m_\varepsilon-l_\varepsilon)^2\right)= 0\,,
\\
 \lim_{\varepsilon \downarrow 0} M_\varepsilon = M =  \Phi_{\xi,\zeta}(\bar t,\bar z,\bar p,\bar m)
\,.
\end{cases}
\end{align}
Moreover, it follows from the boundaries assumptions on $V$ and $U$ that $(\bar{t},\bar{z},\bar p,\bar m)\in(t_i,t_{i+1})\times\R^{d}\times(0,\infty)^{n-i}\times\R^+_*$. As a consequence we assume that, up to a subsequence, $(t_{\varepsilon},z_{\varepsilon},r_{\varepsilon},p_{\varepsilon},q_{\varepsilon},m_{\varepsilon},l_{\varepsilon})\in
(t_i,t_{i+1})\times\R^{d}\times\R^{d}\times[0,\infty)^{(n-i)}\times[0,\infty)^{(n-i)}\times\R^+\times\R^+$.\\
\textbf{Step 2.}
Using Ishii's Lemma (see \cite[Theorem 8.3]{CIL92}) we obtain the existence of real coefficients $\tilde c_{1,\varepsilon},\tilde c_{2,\varepsilon}$, two vectors $\tilde \Delta_{1,\varepsilon},\tilde \Delta_{2,\varepsilon}$ and two symmetric matrices $\tilde{\mathcal{X}}_{\varepsilon}$ and $\tilde{\mathcal{Y}}_{\varepsilon}$ being s.t.
$
(\tilde c_{1,\varepsilon},\tilde \Delta_{1,\varepsilon},\tilde{\mathcal{X}}_{\varepsilon})
\in\bar{\mathcal{J}}^+_{\bar{\mathcal{O}}}({U}_{\xi}(t_{\varepsilon},z_{\varepsilon},p_{\varepsilon},m_{\varepsilon})-\zeta\text{e}^{-\rho t_{\varepsilon}}(1+|z_{\varepsilon}|^2))$ and $(\tilde c_{2,\varepsilon},\tilde\Delta_{2,\varepsilon},\tilde{\mathcal{Y}}_{\varepsilon})
\in\bar{\mathcal{J}}^-_{\bar{\mathcal{O}}}({V}_{\xi}(t_{\varepsilon},r_{\varepsilon},q_{\varepsilon},l_{\varepsilon})+\zeta\text{e}^{-\rho t_{\varepsilon}}(1+|r_{\varepsilon}|^2)),
$
with $\bar{\mathcal{O}}\defeq(t_i,t_{i+1})\times\R^{d}\times(0,\infty)^{n-i}\times \R^+_*$ and $\bar{\mathcal{J}}^+$ (resp. $\bar{\mathcal{J}}^-$) the limiting second-order super-jet (resp. sub-jet) of $U_\xi$ (resp. $V_\xi$) at $(t_{\varepsilon},z_{\varepsilon},p_{\varepsilon},m_{\varepsilon}) \in\bar{\mathcal{O}}$ (resp. $(t_{\varepsilon},r_{\varepsilon},q_{\varepsilon},l_{\varepsilon}) \in\bar{\mathcal{O}}$) and
where
\begin{gather}
\tilde c_{1,\varepsilon}-\tilde c_{2,\varepsilon}= 0,\quad
 \tilde\Delta_{1,\varepsilon}
  =\left(\begin{matrix}\frac{1}\varepsilon (z_{\varepsilon}-r_{\varepsilon})\\\frac{1}\varepsilon (p_{\varepsilon}-q_{\varepsilon})\\\frac{1}\varepsilon (m_{\varepsilon}-l_{\varepsilon})\end{matrix} \right)
\mbox{ and }
\tilde\Delta_{2,\varepsilon}
=\left(\begin{matrix}\frac{1}\varepsilon (z_{\varepsilon}-r_{\varepsilon})\\\frac{1}\varepsilon (p_{\varepsilon}-q_{\varepsilon})\\\frac{1}\varepsilon (m_{\varepsilon}-l_{\varepsilon}) \end{matrix} \right),\quad
 \left(\begin{matrix} \tilde{\cX}_{\varepsilon}&0\\0&-\tilde{\cY}_{\varepsilon} \end{matrix}\right)\leq
\frac{3}{\varepsilon}\left(\begin{matrix} I &-I\\-I& I \end{matrix}\right)
\,.\nonumber
\end{gather}
Hence appealing to  \cite[Remark 2.7 (ii)]{CIL92}, one has the existence of $
( c_{1,\varepsilon}, \Delta_{1,\varepsilon},{\mathcal{X}}_{\varepsilon})
\in\bar{\mathcal{J}}^+_{\bar{\mathcal{O}}}{U}_{\xi}(t_{\varepsilon},z_{\varepsilon},p_{\varepsilon},m_{\varepsilon})$ and $( c_{2,\varepsilon},\Delta_{2,\varepsilon},{\mathcal{Y}}_{\varepsilon})
\in\bar{\mathcal{J}}^-_{\bar{\mathcal{O}}}{V}_{\xi}(t_{\varepsilon},r_{\varepsilon},q_{\varepsilon},l_{\varepsilon}),
$ s.t.
\begin{subequations}\label{eq: Der}
\begin{align}
& c_{1,\varepsilon}=\tilde c_{1,\varepsilon}-\zeta\rho\text{e}^{-\rho t_{\varepsilon}}(1+|z_{\varepsilon}|^{2})\mbox{ and }c_{2,\varepsilon}  = \tilde c_{2,\varepsilon}+\zeta\rho\text{e}^{-\rho t_{\varepsilon}}(1+|r_{\varepsilon}|^{2})\,,\label{eq: Der b1}\\
& \Delta_{1,\varepsilon}=\tilde \Delta_{1,\varepsilon}+
2\zeta\text{e}^{-\rho t_{\varepsilon}} ( z_{\varepsilon},0,0)^\top
\mbox{ and }
\Delta_{2,\varepsilon}=\tilde \Delta_{2,\varepsilon}-
2\zeta\text{e}^{-\rho t_{\varepsilon}}  (r_{\varepsilon},0,0)^\top\,,\label{eq: Der f}
\\
& \left(\begin{matrix} \cX_{\varepsilon}&0\\0&-\cY_{\varepsilon} \end{matrix}\right)=\left(\begin{matrix} \tilde{\cX}_{\varepsilon}&0\\0&-\tilde{\cY}_{\varepsilon} \end{matrix}\right)+2\zeta\text{e}^{-\rho t_{\varepsilon}}\left(\begin{matrix} \bar I &\textbf{0}\\\textbf{0}&\bar I \end{matrix}\right)\leq
\frac{3}{\varepsilon}\left(\begin{matrix} I &-I\\-I& I \end{matrix}\right)+2\zeta\text{e}^{-\rho t_{\varepsilon}}\left(\begin{matrix} \bar I &\textbf{0}\\\textbf{0}&\bar I \end{matrix}\right)
\,,
\end{align}
\end{subequations}
with $I\equiv I_{d+n-i+1}$ and $\bar I$ defined in Lemma \ref{le ModulusCont0}.
Thus, Lemma \ref{le ModulusCont0} and \eqref{eq: Der} imply
\begin{align}
&\sup_{\substack{{(u,b)\in U\times\mathcal{S}_{d(n-i)+d}}}}\mathrm H_{\kappa,\lambda}^{u,b} \left(t_{\varepsilon},r_{\varepsilon},q_{\varepsilon},l_{\varepsilon},
\Delta_{2,\varepsilon},\mathcal{Y}_{\varepsilon},c_{2,\varepsilon}\right)- \sup_{\substack{{(u,b)\in U\times\mathcal{S}_{d(n-i)+d}}}}\mathrm H_{\kappa,\lambda}^{u,b} \left(t_{\varepsilon},z_{\varepsilon},p_{\varepsilon},m_{\varepsilon},\Delta_{1,\varepsilon},\mathcal{X}_{\varepsilon},c_{1,\varepsilon}\right)\nonumber\\
&\quad\quad\quad\quad\quad\quad\quad\quad\quad\quad\quad\quad\quad\quad\le\frac{C}{\varepsilon}\left(|z_{\varepsilon}-r_{\varepsilon}|^2
+|p_{\varepsilon}-q_{\varepsilon}|^2+(m_{\varepsilon}-l_{\varepsilon})^2\right)\label{eq CtradictoryModCont}
\,,
\end{align}
for some $C>0$. Sending  $\varepsilon$ to zero and using \eqref{eq cv}, the last inequality is non-positive.\\
\textbf{Step 3.} We also know from the definition of ${U}_{\xi}$ and ${V}_{\xi}$ that they are respectively sub-/super-solution of \eqref{viscsolBETbis}. As a result, appealing to Lemma \ref{le StrictSupersol} we obtain
\begin{align*}
&\sup_{\substack{{(u,b)\in U\times\mathcal{S}_{d(n-i)+d}}}}\mathrm H_{\kappa,\lambda}^{u,b} \left(t_{\varepsilon},r_{\varepsilon},q_{\varepsilon},l_{\varepsilon},
\Delta_{2,\varepsilon},\mathcal{Y}_{\varepsilon},c_{2,\varepsilon}\right)\nonumber\\&\quad - \sup_{\substack{{(u,b)\in U\times\mathcal{S}_{d(n-i)+d}}}}\mathrm H_{\kappa,\lambda}^{u,b} \left(t_{\varepsilon},z_{\varepsilon},p_{\varepsilon},m_{\varepsilon}, \Delta_{1,\varepsilon},\mathcal{X}_{\varepsilon},c_{1,\varepsilon}\right)
\ge \xi  \frac{1}{8}>0 \,,
\end{align*}
contradicting \eqref{eq CtradictoryModCont}. Hence $(U_{\xi}-V_{\xi})\leq0$ for all $\xi>0$ on $(t_i,t_{i+1})\times\R^{d}\times[0,\infty)^{n-i}\times\R^+$.
\qed

\bibliography{BibGB}
\bibliographystyle{spmpsci_unsrt}
\end{document}